\documentclass[12pt,reqno]{amsart}

\usepackage{amsthm}
\usepackage{amssymb}
\usepackage{latexsym}
\usepackage{multicol}
\usepackage{verbatim,enumerate}
\usepackage{accents}
\usepackage[usenames]{color}
\usepackage{youngtab}
\usepackage[all, poly]{xy}  
\usepackage[utf8]{inputenc}
\usepackage{hyperref}
\usepackage{amsmath, amscd}
\usepackage{soul}
\usepackage{tikz} 
\usetikzlibrary{matrix,arrows,decorations.pathmorphing}
\theoremstyle{definition}

\newtheorem*{defn}{Definition}
\newtheorem*{rem}{Remark}
\newtheorem*{prop}{Proposition}
\newtheorem*{thm}{Theorem}
\newtheorem*{lem}{Lemma}
\newtheorem*{cor}{Corollary}
\newtheorem*{conj}{Conjecture}

\newcommand\bomega{{\mbox{\boldmath $\omega$}}}
\newcommand\bomegas{{\mbox{{\boldmath {\small$\omega$}}}}}

\newcommand{\lie}[1]{\mathfrak{#1}}
\newcommand{\wh}[1]{\widehat{#1}}

\newcommand{\weight}[1]{{\rm{wt}}{#1}}
\newcommand\ba{\bold A}
\newcommand\bc{\mathbb C}
\newcommand\bn{\mathbb N}
\newcommand\bz{\mathbb Z}

\newcommand{\ev}{\operatorname{ev}}
\newcommand{\Hom}{\operatorname{Hom}}

\newcommand{\ch}{\operatorname{ch}}

\newcommand{\qbinom}[2]{\genfrac[]{0pt}0{#1}{#2}}
\newcommand{\floor}[1]{\left\lfloor #1 \right\rfloor}

\def\a{\alpha}

\newcommand{\cal}[1]{\mathcal{#1}}

\def\gr{\operatorname{gr}}
\def\loc{\operatorname{loc}}

\advance\textwidth by 1.2in  \advance\oddsidemargin by -.6in \advance\evensidemargin by -.6in

\newenvironment{pf}{\proof}{\endproof}
\newcounter{cnt}
\newenvironment{enumerit}{\begin{list}{{\hfill\rm(\roman{cnt})\hfill}}{%
\settowidth{\labelwidth}{{\rm(iv)}}\leftmargin=\labelwidth%
\advance\leftmargin by \labelsep\rightmargin=0pt\usecounter{cnt}}}{\end{list}} \makeatletter
\def\mydggeometry{\makeatletter\dg@YGRID=1\dg@XGRID=20\unitlength=0.003pt\makeatother}
\makeatother \theoremstyle{remark}

\numberwithin{equation}{section}
\makeatletter
\def\section{\def\@secnumfont{\mdseries}\@startsection{section}{1}%
  \z@{.7\linespacing\@plus\linespacing}{.5\linespacing}%
  {\normalfont\scshape\centering}}
\def\subsection{\def\@secnumfont{\bfseries}\@startsection{subsection}{2}%
  {\parindent}{.5\linespacing\@plus.7\linespacing}{-.5em}%
  {\normalfont\bfseries}}
\def\subsubsection{\def\@subsecnumfont{\bfseries}\@startsection{subsubsection}{3}%
  {\parindent}{.5\linespacing\@plus.7\linespacing}{-.5em}%
  {\normalfont\bfseries}}
  \makeatother
\newcommand\bu{\bold U}

\begin{document}

\title[]{Quantum affine algebras, graded limits and flags}
\author{Matheus Brito}
\address{Departamento de Matematica, UFPR, Curitiba - PR - Brazil, 81530-015}
\email{mbrito@ufpr.br}
\thanks{}

\author{Vyjayanthi Chari}
\address{Department of Mathematics, University of California, Riverside, 900 University Ave., Riverside, CA 92521, USA}
\email{chari@math.ucr.edu}
\thanks{V.C. was partially supported by DMS-1719357, the Simons Fellows Program and by the  Infosys Visiting Chair position at the Indian Institute of Science.}

\author{Deniz Kus}
\address{University of Bochum, Faculty of Mathematics, Universitätsstr. 150, 44801 Bochum, Germany}
\email{deniz.kus@rub.de}
\thanks{D.K. was partially funded by the Deutsche Forschungsgemeinschaft (DFG, German Research Foundation)-- grant 446246717.}

\author{R. Venkatesh}
\address{Department of Mathematics, Indian Institute of Science, Bangalore 560012, India}
\email{rvenkat@iisc.ac.in}
\thanks{R.V. was partially supported by the Infosys Young Investigator Award grant.}
\subjclass[2010]{}

\begin{abstract}
    In this survey, we  review some of the  recent connections between  the representation theory of (untwisted) quantum affine algebras and the representation theory of current algebras. We mainly focus on the finite-dimensional representations of these algebras. This connection arises via the notion of the   graded and classical limit of finite-dimensional  representations of  quantum affine algebras.  We explain how this study has led to interesting connections with  Macdonald polynomials and  discuss a BGG-type reciprocity result.  We also discuss the role of Demazure modules in this theory and several recent results on the presentation, structure and combinatorics of Demazure modules.

\end{abstract}

\maketitle

\section{Introduction}
Quantized enveloping   algebras  were introduced independently by Drinfeld (1985) and Jimbo (1986) in the context of integrable systems and  solvable lattice models and give a systematic way to construct solutions to the quantum Yang-Baxter equation.   The quantized algebra associated to an affine Lie algebra is called a quantum affine algebra. The representation theory of these  has been intensively studied for nearly thirty-five years since its introduction. It has connections with many research areas of mathematics and physics, for example, statistical mechanics, cluster algebras, dynamical systems, the geometry of quiver varieties, Macdonald polynomials to name a few. In this survey, we mainly focus on the category of finite-dimensional representations $\mathcal{F}_q$ of quantum affine algebras and their connections to graded representations of current algebras. The fact that this  category is not semi-simple  gives a very rich structure and has many interesting consequences. The category is studied via the Drinfeld realization of quantum affine algebras and irreducible objects are parametrized in terms of Drinfeld polynomials. The classical version of $\mathcal{F}_q$ was studied previously in \cite{Cha86}, \cite{CP91}, and the irreducible finite-dimensional representations of the affine algebra and the loop algebra were classified in those papers.

However, we still  have  limited information 
on the structure of finite-dimensional representations of quantum affine algebras except for a few special cases. For example, we do not even know the dimension formulas in general. One  way to study these representations is to go from quantum level to classical level by forming the classical limit, see for instance \cite{CP01} for a necessary and sufficient condition for the existence of the classical limit. The classical limit (when it exists) is a finite-dimensional module for the corresponding affine Lie algebra. By restricting and suitably twisting this classical limit, we obtain the graded limit which is a graded representation of the corresponding current algebra, see Section \ref{classgrad} for more details. 
Most of the time we get a reducible indecomposable representation of affine Lie algebra (or current algbera) on passing to the classical (graded) limit. A similar phenomenon is observed in modular representation theory: an irreducible finite-dimensional representation in characteristic zero becomes reducible on passing to characteristic $p$. 
Many interesting families of representations from $\mathcal{F}_q$ admit this graded limit, for instance, the local Weyl modules, Kirillov-Reshetikhin modules, minimal affinizations, and some of the  prime representations coming from the work of Hernandez and Leclerc on monoidal categorification. 

In \cite{CP01} the authors introduced  the notion of   local Weyl modules for a  quantum affine algebra. They are given by generators and relations, are highest weight modules in a suitable sense and satisfy a canonical  universal property. In particular, any irreducible module in $\mathcal{F}_q$ is  a quotient of some Weyl module. It was conjectured in \cite{CP01} (and proved there for $\lie{sl}_2$) that any local Weyl module has a tensor product decomposition into fundamental local Weyl modules, see Section \ref{weyltp} for more details. This conjecture stimulated a lot of research on this topic and the general case was established through the work of \cite{CL06, FoL07, Na11}. The work of Kirillov and Reshetikhin \cite{KR87} has a connection with the irreducible representations of quantum affine algebras corresponding to a multiple of a fundamental weight. These modules are referred as Kirillov-Reshetikhin modules in the literature, because they conjectured the classical decomposition of these modules in their paper. The study of Kirillov-Reshetikhin modules has been of immense interest in recent years due to their rich combinatorial structures and several applications to mathematical physics \cite{HKOTY99, He06}. Many important conjectures on the character formulas of these modules and their fusion products were made from physical considerations and they stimulated lots of research, see \cite{FK08, Na17, OS08} and the references therein.

One of the very natural questions that arises from the work of \cite{FR92} is: what is the smallest representation from $\mathcal{F}_q$ that corresponds to a given finite-dimensional irreducible representation of the underlying simple Lie algebra $\lie g$? The second author introduced the notion of a minimal affinization in \cite{Cha95a} with this motivation and it was further studied in \cite{Ch95b,CP95a,CP95}. One introduces a poset for each dominant integral weight, such that each element of the poset determines a family of irreducible representations in $\mathcal{F}_q$. The irreducible modules that correspond to the minimal elements of this poset are minimal affinizations. Kirillov-Reshetikhin modules are the minimal affinizations of multiples of the fundamental weights.
Our final example of graded limits comes from the work of Hernandez and Leclerc  \cite{HL10} on the monoidal categorification of cluster algebras. The authors defined an interesting monoidal subcategory of $\mathcal{F}_q$ in simply-laced type and proved that for $\mathfrak{g}$ of type $A_n$ and $D_4$, it categorifies a cluster algebra of the same type, i.e., its Grothendieck ring admits a cluster algebra structure of the same type as $\mathfrak{g}$. The prime real representations of this subcategory are the cluster variables and these are called the HL-modules. A more detailed discussion of HL-modules can be found in Section~\ref{hlm}.

Even though the study of graded representations of current algebras is mainly motivated by their connection with the representations of quantum affine algebras (via graded limits), they are now of independent interest as they have found many applications in number theory, combinatorics, and mathematical physics. They have connections with mock theta functions, cone theta functions \cite{BCSV15, BCSW16, BKu17}, symmetric Macdonald polynomials \cite{BCSW21, CI15, I03}, the $X = M$ conjecture \cite{AK07, FK08, Na17}, and Schur positivity \cite{FHD14, V13} etc.
One of the very important families of graded representations of current algebras comes from $\lie g$-stable Demazure modules. A Demazure module by definition is a module of the Borel subalgebra of the affine Lie algebra. If it is $\lie g$-stable, then it naturally becomes a module of the maximal parabolic subalgebra which contains the current algebra. By restriction, we get a graded module of the current algebra. These modules are parametrized by pairs consisting of a positive integer and a dominant weight, and given such a pair $(\ell, \lambda)$, the corresponding $\lie g$-stable Demazure module of the current algebra is denoted by $D(\ell, \lambda)$. These modules include all well-known families of graded representations of current algebras. For example, any local Weyl module of the current algebra is isomorphic to a level one Demazure module $D(1, \lambda)$ when $\lie g$ is simply-laced.

The limit of a tensor product of quantum affine algebra modules is not necessarily isomorphic to the tensor product of their classical limits. So, we need to replace the tensor product with something else in order to study the limit of a tensor product of quantum affine algebra modules. Examples suggest that the fusion product introduced by Feigin and Loktev \cite{FL99} is the correct notion that should replace the tensor product. It is a very important and seemingly very hard problem to understand the fusion products of $\lie g$-stable Demazure modules of various levels. One would like to find the generators and relations and the graded character of these modules, but very limited cases are known \cite{BK20, CV13, F15a,  Na17}.

The survey is organized as follows. We begin by stating the foundational results, including the definition of local Weyl modules and the classification of irreducible modules in Section 2. In Section 3, we discuss various well-studied families of finite-dimensional representations of quantum affine algebras and review the literature on the presentation of these modules, their classical limit, and the closely related graded limits. Later we  move on to the study of graded finite-dimensional representations of current algebras. We relate the local Weyl modules to the $\lie g$-stable Demazure modules and discuss the connection between the characters of the local Weyl modules and Macdonald polynomials. We also discuss the BGG-type reciprocity results and briefly mention some recent developments on tilting modules, generalized Weyl modules, and global Demazure modules. In the end, we collect together some results on Demazure modules.\\\\
{\em Acknowledgements. Part of this paper was written while the second author was visiting the Max Planck Institute, Bonn, in Fall 2021. She gratefully acknowledges the financial support and the excellent working conditions provided by  the institute.}

\section{The simple and untwisted affine Lie algebras}\label{simaff} In this section we collect the notation and some well-known results on the structure and representation theory of affine Lie  algebras.

\subsection{Conventions}\label{conv} We let $\mathbb C$ (resp. $\mathbb Q$, $\mathbb Z$, $\mathbb Z_+$, $\mathbb N$) be the set of complex numbers (resp. rational numbers, integers, non-negative integers, positive integers). We adopt the convention that given two complex vector spaces $V,W$  the corresponding tensor product $V\otimes_{\mathbb C}W$ will be just denoted as $V\otimes W$.  
\\\\ Given an indeterminate $t$ we let $\mathbb C[t]$ (resp. $\mathbb C[t,t^{-1}]$, $\mathbb C(t)$) be the ring of polynomials (resp. Laurent polynomials, rational functions) in the variable $t.$ For $s\in\mathbb{Z}, m,r\in \mathbb Z_+$ with $m\ge r$, set 
$$[s]_t=\frac{t^{s}-t^{-s}}{t-t^{-1}},\ \ [m]_t!=[m]_t[m-1]_t\cdots [1]_t,\ \ {\qbinom{m}{r}}_t=\frac{[m]_t!}{[r]_t![m-r]_t!}.$$
For any complex Lie algebra $\lie a$ we let $\bu(\lie a)$ be the corresponding universal enveloping algebra. Given any commutative associative algebra $A$ over $\mathbb{C}$ define a Lie algebra structure on $\lie a\otimes A $ by $$[x\otimes a, y\otimes b]=[x,y]\otimes ab,\ \ x,y\in\lie a, \ a,b\in A.$$
In the special case when $A$ is $\mathbb C[t]$ or $\mathbb C[t,t^{-1}]$ we set $$\lie a[t]=\lie a\otimes \mathbb C[t],\ \ L(\lie a)=\lie a\otimes \mathbb C[t^{\pm 1}].$$
\subsection{Simple and affine Lie algebras}
\subsubsection{The simple Lie algebra $\lie g$}\label{simples} Let $\lie g$ denote a simple finite-dimensional Lie algebra over $\mathbb C$  and let $\lie h$ be a fixed Cartan subalgebra of $\lie g$ and $R$ the corresponding set of roots. Let $I=\{1,\dots, n\}$ be an index set for the set of simple roots $\{\alpha_1,\dots,\alpha_n\}$ of $R$ and $\{\omega_1,\dots,\omega_n\}$ a set of fundamental weights.  Given $\lambda,\mu\in\lie h^*$ we say that $$\lambda\ge \mu\iff \lambda-\mu\in\sum_{i\in I}\mathbb Z_+\alpha_i.$$
Let $P$, $Q$ (resp. $P^+, Q^+$) be the $\mathbb Z$-span (resp. $\mathbb Z_+$-span) of the fundamental weights and simple roots respectively and let $R^+= R\cap Q^+$. We denote by $\theta\in R^+$ the highest root in $R^+$ and let $(\ , \ )$ be the form on $\lie h^*$ induced by the restriction of the Killing form $\kappa:\lie g\otimes\lie g\to\mathbb C$ of $\lie g$. We assume that it is  normalized so that $(\theta,\theta)=2$  and set $$d_\alpha= 2/(\alpha,\alpha),\ \ d_i= d_{\alpha_i},\ \ \alpha_i^\vee= d_i\alpha_i,\ \ \omega_i^\vee= d_i\omega_i,\ \ a_{i,j}= (\alpha_j, \alpha_i^\vee),\ \ 1\le i,j\le n.$$

Let $W$ be the Weyl group of $\lie g$; recall that it is the subgroup of ${\rm{ Aut}}(\lie h^*)$ generated by the reflections $s_i, i\in I$, defined by:$$s_i(\lambda)=\lambda-(\lambda, \alpha_i^\vee)\alpha_i,\ \ i\in I.$$ 
Fix a Chevalley basis $\{x_\alpha^\pm,\ h_i:\alpha\in R^+, i\in I\}$ of $\lie g$ and set for simplicity $x_i^\pm=x_{\alpha_i}^\pm$. The elements $x_i^\pm, h_i$, $i\in I$ generate $\lie g$ as a Lie algebra. 
Given $\alpha=\sum_{i=1}^nr_i\alpha_i\in R^+$ let $h_\alpha\in\lie h$ be given by $ h_\alpha=d_\alpha\sum_{i=1}^n  \frac{r_i}{d_i}h_i$ and note that the elements $s_{\alpha}, \alpha\in R^+$, defined by $s_{\alpha}(\lambda)=\lambda-\lambda(h_\alpha)\alpha$ are elements of $W$ and we have $s_i=s_{\alpha_i}$.\\\\
Let $\lie n^\pm$ be the subalgebra generated by the elements $\{x_i^\pm: i\in I\}$. Then,  
$$\lie n^\pm=\bigoplus_{\alpha\in R^+}\mathbb Cx_\alpha^\pm,\ \ \ \  \lie b^{\pm}=\lie h\oplus\lie n^\pm,\ \ \lie g= \lie b^\pm\oplus\lie n^\mp.$$ 
We have a corresponding decomposition of $\bu(\lie g)$ as vector spaces $$\bu(\lie g)\cong\bu(\lie n^-)\otimes \bu(\lie h)\otimes \bu(\lie n^+) \cong\bu(\lie n^-)\otimes \bu(\lie b^+).$$

\subsubsection{The affine Lie algebra}\label{affines} The (untwisted) affine Lie algebra $\widehat{\lie g}$ and its Cartan subalgebra $\widehat{\lie h}$ are  defined as follows: $$\widehat{\lie g}= L(\lie g)\oplus\mathbb Cc\oplus\mathbb  Cd,\ \ \ \widehat{\lie h}=\lie h\oplus\mathbb Cc\oplus\mathbb Cd$$ with commutator given by requiring $c$ to be central and $$[x\otimes t^r, y\otimes t^s]=[x,y]\otimes t^{r+s}+r\delta_{r+s,0}\kappa(x,y)c,\ \ [d,x\otimes t^r]=r x\otimes t^r.$$ Here $x,y\in\lie g$, $r,s\in\mathbb Z$ and $\delta_{n,m}$ is the Kronecker delta symbol. Setting $h_0=-h_\theta+c$ we see that the set $\{h_i, d:0\le i\le n\}$ is a basis for $\widehat{\lie h}$. \\\\
Regard an element $\lambda\in\lie h^*$ as an element  of $\widehat{\lie h}^*$ by setting $\lambda(c)=0=\lambda(d)$. Define elements $\delta$, $\alpha_0$ and the affine fundamental weights $\Lambda_i$, $0\le i\le n$,  of $\widehat{\lie h}^*$ by:
$$\delta(d)=1,\ \ \delta(\lie h\oplus \mathbb Cc)=0,\ \ \alpha_0=-\theta+\delta,\ \ \Lambda_0(c)=1,\ \ \Lambda_0(\lie h\oplus\mathbb Cd)=0,$$
$$\Lambda_i(h_j)=\delta_{i,j},\ \ \Lambda_i(d)=0,\ \ i\in I,\ \  j\in\{0,\dots,n\}.$$
The subset $$\widehat R=\{\alpha+r\delta:\alpha\in R\cup\{0\}, r\in\mathbb Z\}\backslash\{0\}\subseteq\widehat{\lie h}^*,$$ is called the set
of affine roots. The set of affine simple roots is $\{\alpha_i:i\in\widehat I\}$ where $\widehat I=\{0,1,\dots ,n\}$. The corresponding set of positive roots is given by:
 $$\widehat R^+=\{\pm\alpha+(r+1)\delta:\alpha\in R^+\cup\{0\}, r\in\mathbb Z_+\}\ \cup R^+.$$   Set
 $$\widehat P=\sum_{i=0}^n \mathbb{Z}\Lambda_i+\mathbb{Z}\delta,\ \ \ \widehat P^+=\sum_{i=0}^n \mathbb{Z}_+\Lambda_i+\mathbb{Z}\delta$$ Let $\widehat Q$ (resp. $\widehat Q^+)$ be the $\mathbb Z$-span (resp. $\mathbb Z_+$-span) of the affine simple roots.
The affine Weyl group $\widehat W$ is the subgroup of  ${\rm{Aut}}(\widehat{\lie h}^*)$ generated by the set $\{s_i: i\in \widehat I\}$ where $$s_i(\lambda)=\lambda-\lambda(h_i)\alpha_i, \ \ i\in\widehat I,\ \ \lambda\in\widehat {\lie h}^{*}.$$
Clearly $W$ is a subgroup of $\widehat W$ and we have an isomorphism $$\widehat W\cong W\ltimes \sum_{i\in I}\mathbb Z \alpha^{\vee}_i.$$ It will also be convenient to introduce the extended affine Weyl group $\widetilde W= W\ltimes \sum_{i\in I}\mathbb Z\omega^{\vee}_i.$ 
\\\\ Setting $x_0^\pm =x_\theta^\mp\otimes t^{\pm 1}$ we observe that $\widehat{\lie g}$ is generated by the elements $\{x_i^\pm, h_i: i\in\widehat I\}\cup\{d\}$.
The root space corresponding to an element $\pm \alpha+s\delta \in \widehat R$ with $\alpha\in R^+$  is  $\mathbb C(x^\pm_\alpha\otimes t^s)$  and to an element $r\delta$ is $(\lie h\otimes t^r)$, $s,r\in\mathbb Z$ and $r\ne 0$.  We shall just denote a non-zero element of the one-dimensional root space corresponding to $\pm\alpha$, $\alpha\in \widehat R^+\setminus\bn\delta$ by $x_\alpha^\pm$ and let $h_\alpha=[x_\alpha^+,x_\alpha^-]$.
The subalgebras $\widehat{\lie n}^\pm$ and $\widehat{\lie b}$ are defined in the obvious way and we have $$\widehat{\lie n}^{\pm}=\lie g\otimes t^{\pm 1}\mathbb C[t]\oplus \lie n^{\pm},\ \ \widehat{\lie b}^+=\widehat{\lie h}\oplus\widehat{\lie n}^+.$$ 
This gives rise to an analogous triangular decomposition $$\bu(\widehat{\lie g})\cong  \bu(\widehat{\lie n}^-)\otimes \bu(\widehat{\lie h})\otimes \bu(\widehat{\lie n}^+)\cong \bu(\widehat{\lie n}^-)\otimes \bu(\widehat{\lie b}^+).$$

\subsubsection{The loop algebra $L(\lie g)$ and the current algebra $\lie g[t]$} 
It is trivial to see that $L(\lie g)\oplus\mathbb Cd$ is the quotient of $\wh{\lie g}$ by the center $\mathbb Cc$. The action of  $d$ obviously induces a $\mathbb Z$-grading on $L(\lie g)$ and the current algebra $\lie g[t]$ is a graded subalgebra of $L(\lie g)$. 
Moreover  $\lie g[t]\oplus\mathbb Cc\oplus\mathbb Cd$  can also be regarded as a maximal parabolic subalgebra of $\widehat{\lie g}$, namely $$\lie g[t]\oplus\mathbb Cc\oplus\mathbb Cd \cong \widehat{\lie b^+}\oplus \lie n^-.$$
We make the grading on $L(\lie g)$ and $\lie g[t]$ explicit for the reader's convenience. For $r\in\mathbb Z$ we declare $\lie g\otimes t^r$ to be the $r$-th graded piece and note that for $r,s\in\mathbb Z$, we have $[\lie g\otimes t^r, \lie g\otimes t^s]=\lie g\otimes t^{r+s}$. This grading induces a grading on the corresponding enveloping algebras as well once we declare a monomial of the form $(a_1\otimes t^{r_1})\cdots(a_p\otimes t^{r_p})$, $a_s\in\lie g$, $r_s\in\mathbb Z$, $1\le s\le p$ to have grade $(r_1+\cdots+r_p)$.
\subsubsection{Ideals in $L(\lie g)$}\label{ideals} The affine Lie algebra is clearly not simple; the center spans a one-dimensional ideal. One can prove using the explicit realization that this and the derived algebra $L(\lie g)\oplus\mathbb Cc$ are the only non trivial proper ideals in $\wh{\lie g}$.
The following result is well-known, a proof can be found for instance in  \cite[Lemma 1]{CFK10}.
\begin{lem} For all $f\in\mathbb C[t,t^{-1}]$ the subspace $\lie g\otimes f\bc[t,t^{-1}]$ is an ideal in $L(\lie g)$. Moreover any ideal in $L(\lie g)$  must be of this form. In particular all ideals are of finite codimension.   Writing $f=(t-a_1)^{r_1}\cdots (t-a_k)^{r_k}$ with $a_r\ne a_s$ for $1\le r\ne s\le k$ we see that we have an isomorphism of Lie algebras $$\frac{L(\lie g)}{\lie g\otimes f\bc[t,t^{-1}]}\cong \lie g\otimes \frac{\mathbb C[t,t^{-1}]}{(f)}\cong \bigoplus_{s=1}^k \left(\lie g\otimes \frac{\bc[t,t^{-1}]}{(t-a_s)^{r_s}}\right).$$
\hfill\qedsymbol
\end{lem}
We shall sometimes refer to the finite-dimensional quotient of $L(\lie g)$ defined by $f\in\mathbb C[t,t^{-1}]$ as the  truncation of $L(\lie g)$ at $f.$

\subsection{Representations of simple and affine Lie algebras}
\subsubsection{Finite-dimensional representations of $\lie g$}\label{dimfing} We say that a $\lie g$-module $V$ is a weight module if, 
$$V=\bigoplus_{\mu\in\lie h^*} V_\mu,\ \ V_\mu=\{v\in V: hv=\mu(h)v,\ \forall h\in\lie h\}$$ and we let ${\rm {wt}}(V)=\{\mu\in\lie h^*: V_\mu\ne 0\}.$ If  ${\rm wt}(V) \subset P$ and $\dim V_\mu<\infty$ for all $\mu\in P$,  we let ${\rm{ch}}(V)$ be the formal sum $\sum_{\mu\in P}(\dim V_\mu)\  e_\mu$ where $e_\mu$ varies over a basis of the group ring $\mathbb Z[P]$. \\\\
Given $\lambda\in P^+$  let $V(\lambda)$ be the $\lie g$-module generated by an element $v_\lambda$ with defining relations:
$$h_iv_\lambda=\lambda(h_i)v_\lambda,\ \  x_i^+ v_\lambda=0, \ \ (x_i^-)^{\lambda(h_i)+1} v_\lambda=0,\ \ 1\le i\le n.$$
It is well-known \cite{Hu80} that $V(\lambda)$ is a weight module with ${\rm {wt}}(V(\lambda))\subseteq \lambda- Q^+$ and that it is also an irreducible and  finite-dimensional $\lie g$-module.  The set ${\rm {wt}}(V(\lambda))$ is $W$-invariant and $\dim V(\lambda)_\mu=\dim V(\lambda)_{w\mu}$ for all $w\in W$. \\\\
Any finite-dimensional $\lie g$-module is isomorphic to a direct sum of copies of  $V(\lambda)$, $\lambda\in P^+$. 
\subsubsection{ Integrable and positive level representations of $\widehat{\lie g}$}\label{levelr}   The notion of a weight module and its character for $\widehat{\lie g}$ are defined as for  simple Lie algebras with $\lie h$ replaced by $\widehat{\lie h}$.
We say that  a  $\widehat{\lie g}$-module  $V$ is integrable if it is a weight module and the elements $x_i^\pm $, $i\in\widehat I$, act locally nilpotently.
We say that $V$   is of  level  $r\in\mathbb Z$ if   $cv=rv$ for all $v\in V$; if $r>0$ (resp. $r<0$) then we say that $V$ is of positive level (resp. negative level).   \\\\
Given $\lambda\in \widehat P^+$  let $V(\lambda)$ be the $\widehat{\lie g}$-module generated by an element $v_\lambda$ with defining relations:
$$h_iv_\lambda=\lambda(h_i)v_\lambda,\ \  x_i^+ v_\lambda=0, \ \ (x_i^-)^{\lambda(h_i)+1} v_\lambda=0,\ \ 0\le i\le n.$$
Again it is well-known  \cite{K90} that $V(\lambda)$ is an integrable irreducible module with positive level $\lambda(c)$ and  $$\dim V(\lambda)_\mu\ne 0\implies\mu\in\lambda-\widehat Q^+,\ \ \dim V(\lambda)_\mu<\infty,\ \ \mu\in\widehat P,$$ $$\dim V(\lambda)_\mu=\dim V(\lambda)_{w\mu},\ \ w\in\widehat W,\ \  \mu\in\widehat P.$$ Notice that the preceding properties show immediately that $V(\lambda)$ is infinite-dimensional as long as $\lambda\notin \mathbb{Z}\delta$. In the case when $V$ is irreducible the following was proved in \cite{Cha86}. The complete reducibility statement was proved in \cite{Es03}.
\begin{thm} Any  positive level integrable $\wh{\lie g}$-module $V$ with the property that $\dim V_\mu<\infty$ for all $\mu\in \widehat P$ is isomorphic to a direct sum of modules of the form $V(\lambda)$, $\lambda\in\widehat P^+$.\hfill\qedsymbol\end{thm}
There is a completely similar theory for negative level modules.
\subsubsection{Demazure modules} \label{demazure} Given  $\lambda\in P^+$ and $w\in W$ (resp. $\lambda\in\widehat P^+, w\in\widehat W$)  the Demazure module $V_w(\lambda)$ is the $\lie b$-submodule (resp. $\wh{\lie b}$-submodule) of $V(\lambda)$ generated by the one dimensional subspace $V(\lambda)_{w\lambda}$. The Demazure modules are always finite-dimensional; in the case of $\lie g$ this statement is trivial while for $\widehat{\lie g}$ the statement follows from the fact that ${\rm wt}( V(\lambda))\subseteq\lambda-\widehat Q^+$.
For $\lambda\in P^+$ and $w\in W$,  (resp.  $\lambda\in \widehat{P}^+$, $w\in \widehat{W}$) it is a result from \cite{Jos85, Ma88} that as a $\bu(\lie b)$-module  (resp. $\bu(\widehat{\lie b})$-module)  $V_w(\lambda)$   is generated by $v_{w\lambda}$ with relations: $hv_{w\lambda}= (w\lambda)(h) v_{w\lambda},$ for all $h\in\lie h$ (resp. for all $h\in\widehat{\lie h}$) and $$ (x_\alpha^+ )^{p+1} v_{w\lambda}=0,\ \ p\ge \max\{0,-w\lambda(h_\alpha)\},\ \ \alpha\in R^+$$
$$({\rm{resp. }}  (\lie h\otimes t^{r+1})v_{w\lambda}=0,\ \ r\in\mathbb Z_+,\ (x^+_\alpha)^{p+1} v_{w\lambda}=0,\ \ p\ge \max\{0,-w\lambda(h_\alpha)\},\ \  \alpha\in \widehat R^+\setminus\mathbb N\delta).$$
These relations were simplified in \cite{CV13} and \cite{KV21} and we refer to Theorem~\ref{mainthmgstable} for a precise statement.

\begin{lem}\label{stable} Suppose that $\lambda\in P^+$,  $w\in W$, and $i\in I$ are such that $(w\lambda)(h_i)\leq 0$. Then $x_i^-V(\lambda)_{w\lambda}=0.$ An analogous statement holds for $V(\lambda)$ with $\lambda\in\widehat P^+ $.\end{lem}
\begin{pf} If $w\lambda-\alpha_i\in {\rm wt}(V(\lambda))$, then $\lambda-w^{-1}\alpha_i\in {\rm wt}(V(\lambda))$. On the other hand our assumptions force $(w\lambda)(h_i)=0$ or $w^{-1}\alpha_i\in R^-$ where the latter condition ends in a contradiction to ${\rm wt}(V(\lambda))\subseteq\lambda-Q^+$. Hence $(w\lambda)(h_i)=0$ and $x_i^-V(\lambda)_{w\lambda}=0$ follows.
\end{pf}
As a consequence of this lemma we see immediately that 
$$\lambda\in\widehat P^+,\ \ w\in \widehat W, \ \ (w\lambda)(h_i)\in -\mathbb Z_+, \ \forall\ i\in I\implies\ \ \lie g[t] V_w(\lambda)\subseteq V_w(\lambda).$$ We call these $\widehat{\lie b}$-submodules the $\lie g$-stable Demazure modules. 
 
\subsubsection{Level zero  modules for $\widehat{\lie g}$ and finite-dimensional modules for $L(\lie g)$ and $\lie g[t]$}\label{levelzero} A level zero module for $\widehat{\lie g}$ is one on which the center acts trivially, in particular it can be regarded  as a module for $L(\lie g)\oplus\mathbb Cd$. More generally given any representation $V$ of $\lie g$ one can define a $L(\lie g)\oplus\mathbb Cd$-module structure on $L(V)=V\otimes\mathbb C[t, t^{-1}]$ by:
$$(x\otimes t^r)(v\otimes t^s)=xv\otimes t^{r+s},\ \ d(v\otimes t^r)=rv\otimes t^r,\ \ r,s\in\mathbb Z,\ x\in\lie g.$$
The only finite-dimensional representations of $\widehat{\lie g}$  are those on which $L(\lie g)\oplus\mathbb C c$ acts trivially. 
We give a proof of this fact for the reader's convenience. Note that by working with the Jordan--Holder series it suffices to prove this for irreducible finite-dimensional representations.
Thus, let  $V$ be  a finite-dimensional irreducible representation of $\wh{\lie g}$.
Then it is easily seen that there exists a vector $0\ne v\in V$ such that the following hold:
$$(x_\alpha^+\otimes t^r)v=0,\ \ (h_i\otimes t^{s})v=a_{i,s}v,\ \ \alpha\in R^+, \ \ i\in I,  r\in \mathbb Z,\ \ s\in\mathbb Z_+, \ \ cv=\ell v, \ \ dv= av,$$
where $a, a_{i,s}\in\mathbb C$ for $s>0$ and $a_{i,0}\in\mathbb Z_+$ and $\ell\in\mathbb Z $. Then $$[d, h_i\otimes t^s]v=s(h_i\otimes t^s)v\implies 0=sa_{i,s}v\implies a_{i,s}=0,\ \ s>0.$$
Working with the $\lie{sl}_2$-triple  $\{x_i^+\otimes t^{-s}, x_i^-\otimes t^{s}, h_i-s\kappa(x_i^+, x_i^-)c\}$ we see that we must have $a_{i,0}-s\ell\kappa(x_i^+,x_i^-)\ge 0$ for all $s\ge 0$ and this forces $\ell=0$. Hence $V$ must be a finite-dimensional representation of $L(\lie g)\oplus\mathbb Cd$. Suppose that $a_{i,0}>0$ for some $i\in I$. Then working again with the triple $\{x_i^+\otimes t^{-s}, x_i^-\otimes t^{s}, h_i-s\kappa(x_i^+,x_i^-)\}$ we see that $(x_i^-\otimes t^s)v\ne 0$ for all $s>0$. Since these elements have $d$-eigenvalues $a+s$ they must be linearly independent which is a contradiction. It follows that $a_{i,0}=0$ for all $i\in I$ and then it is  easily seen that $V$ must be the trivial $L(\lie g)\oplus\mathbb Cc$-module.
\\\\
Consider however, the commutator subalgebra $L(\lie g)\oplus\mathbb Cc$ of $\widehat{\lie g}$. The preceding arguments show that the center must act trivially on any finite-dimensional representation. Hence it suffices to study finite-dimensional representations of $L(\lie g)$. To construct examples we  introduce for each $a\in\mathbb C^\times$ the evaluation homomorphism $\ev_a: L(\lie g)\to \lie g$ which sends $x\otimes t^r\to a^rx$ for $x\in\lie g$ and $r\in\mathbb Z$. Given a representation $V$ of $\lie g$ let $\ev_a V$ be the pull-back $L(\lie g)$-module. The following was proved in \cite{CP86,Cha86}. 
\begin{prop}\hfill
\begin{enumerate}
     \item[(i)]
Any irreducible finite-dimensional representation of $L(\lie g)$ is isomorphic to a tensor product of the form $\ev_{a_1}V(\lambda_1)\otimes\cdots\otimes \ev_{a_k} V(\lambda_k)$ for some $k\ge 1$, pairwise distinct elements $a_1,\dots, a_k\in\mathbb C^\times$, and elements $\lambda_1,\dots,\lambda_k\in P^+$. Moreover any tensor product of irreducible finite-dimensional representations is either irreducible or completely reducible.
\item[(ii)] Suppose that $V$ is an irreducible finite-dimensional module of $L(\lie g)$. Then $L(V)$ is a direct sum of irreducible modules for $L(\lie g)\oplus\mathbb Cd$. Any level zero integrable irreducible  module for $\widehat{\lie g}$ with finite-dimensional weight spaces is obtained as a direct summand of $L(V)$. 

\item[(iii)] A similar result holds for finite-dimensional modules of $\lie g[t]$ once we also allow the module $\ev_0 V(\lambda)$. \hfill\qedsymbol
\end{enumerate}
\end{prop}
Level zero modules for $\wh{\lie g}$ are not completely reducible. The simplest example is the adjoint representation where the center is a proper submodule which does not have a complement.
Finite-dimensional modules for $L(\lie g)$ and $\lie g[t]$ are also not completely reducible. For instance the $\lie g$-stable Demazure modules  are usually reducible and indecomposable. A simple exercise, left to the reader,  is to verify this in the case when  $\widehat{\lie g}=\widehat{\lie{sl}_2}$, $\lambda=\Lambda_0$ and $w=s_1s_0$.\\\\
We shall frequently be interested in $\mathbb Z$-graded modules for $\lie g[t]$. These are $\mathbb Z$-graded vector spaces $V=\bigoplus_{m\in\mathbb Z} V[m]$  which admit an action of $\lie g[t]$ satisfying $(\lie g\otimes t^r)V[m]\subseteq V[m+r]$ for all $m\in\mathbb Z$ and $r\in\mathbb{Z}_+$. In particular each graded piece $V[m]$ is a module for $\lie g$. If $\dim V[m]<\infty$ for all $m\in\mathbb{Z}$ we define the graded character as the formal sum
$$\ch_{\gr}V=\sum_{m\in\mathbb Z} \ch_{\lie g} V[m] q^m.$$ 

The module $\ev_a V(\lambda)$ for $\lie g[t]$ is graded if and only if $a=0$. The $\lie g$-stable Demazure modules are also graded where the grading is given by the action of $d$; namely $V_w(\lambda)[r]=\{v\in V_w(\lambda): dv=rv\}.$ \\\\
Given a $\mathbb Z$-graded vector space $V$ and an integer $s\in\mathbb Z$ let $\tau_sV$ be the grade shifted vector space obtained by declaring $(\tau_sV)[r]=V[r-s].$

\subsubsection{The monoid $\mathcal P^+$ and an alternate parametrization of finite-dimensional modules}\label{drinpar} Let $\mathcal P^+$ be the free abelian multiplicative monoid generated by elements $\bomega_{i,a}$, $i\in I$, $a\in\mathbb C^\times$ and let $\mathcal P$ be the corresponding group generated by these elements. %Let $\weight:\ \mathcal P\to P$ be defined by extending the assignment $\weight(\bomega_{i,a})=\omega_i$ to a homomorphism of abelian groups.

For $\lambda\in P^+$ let $\bomega_{\lambda,a}=\prod_{i=1}^n\bomega_{i,a}^{\lambda(h_i)}.$ Clearly any element of $\cal P^+$ can be written uniquely as a product $\bomega_{\lambda_1,a_1}\cdots\bomega_{\lambda_k,a_k}$ for some multisubset $\{\lambda_1,\dots,\lambda_k\}\subseteq P^+$ and distinct elements $a_s\in\mathbb C^\times$, $1\le s\le k$.  Then part (i) of Proposition \ref{levelzero} can be reformulated as follows.
\begin{prop}
There exists a bijective correspondence between  $\mathcal P^+$ and isomorphism classes of finite-dimensional irreducible representations of $L(\lie g)$  given by  $$\bomega= \bomega_{\lambda_1,a_1}\cdots\bomega_{\lambda_k,a_k}\longrightarrow [V(\bomega)]= [\ev_{a_1}V(\lambda_1)\otimes\cdots\otimes~   \ev_{a_k} V(\lambda_k)].$$
Moreover if $\bomega,\bomega'\in\cal P^+$ then
  $V(\bomega)\otimes V(\bomega')$
is completely reducible and has $V(\bomega\bomega')$ as a summand.
\hfill\qedsymbol
\end{prop}
\subsubsection{Annihilating ideals}\label{trunc} Suppose that $V$ is a finite-dimensional representation of $L(\lie g)$. Then the discussion in Section \ref{ideals} shows that there exists $f\in\mathbb C[t,t^{-1}]$ such that $$\{a\in L(\lie g): av=0\ \ {\rm{for\ all}}\ \ v\in V\}=\lie g\otimes f\mathbb C[t,t^{-1}].$$
In particular $V$ becomes a module for the truncated  Lie algebra $\lie g\otimes \mathbb C[t,t^{-1}]/(f)$.  In the case when $V$ is irreducible the discussion so far proves that $f=(t-a_1)\cdots (t-a_k)$ for some distinct element $a_1,\dots, a_k\in \mathbb C^\times$.  \\\\
Suppose that $V_1$ and $V_2$ are modules for the truncation of $L(\lie g)$ at $f_1$ and $f_2$ respectively and let $f$ be the least common multiple of the pair. Then $V_1\otimes V_2$ is a module for the truncation at $f$.
\section{The quantized simple and affine enveloping algebras} \label{quantum}
\subsection{Definitions and the Hopf algebra structure}
\subsubsection{The Drinfeld-Jimbo presentation} Let $q$ be an indeterminate and set $q_i = q^{d_i}$, $i\in I$, and $q_0 = q$. The quantized enveloping algebra $\bu_q(\widehat{\lie g})$ (also called the quantum affine algebra) is the associative algebra over $\mathbb C(q)$  generated by elements $X_i^\pm, K_i^{\pm 1}$, $D^{\pm 1}$, $i\in\widehat I$ and relations:
\begin{gather*} 
K_iK_i^{-1} = 1, \ \ \ DD^{-1}= 1, \ \ \  K_iK_j = K_jK_i, \ \ \ DK_i=K_iD, \ \  \ i,j\in \widehat I\\ 
K_iX_j^{\pm}K_i^{-1}=q_i^{\pm a_{ij}} X_j^\pm,\ \ DX_j^{\pm}D^{-1}=q^{\pm \delta_{0,j}} X_j^\pm, \ \ \ \  
[X_i^+,X_j^-] = \delta_{i,j}\frac{K_i-K_i^{-1}}{q_i-  q_i^{-1}}, \ \ \ i,j\in \widehat I,\\ 
\sum_{m=0}^{1-a_{i,j}} (-1)^m {\qbinom{1-a_{i,j}}{m}}_{q_i} (X_i^\pm)^{1-a_{i,j}-m}X_j^\pm (X_i^\pm)^{m}=0, \ \ \ i,j\in \widehat I, \ i\neq j.
 \end{gather*}
The Hopf structure on this algebra is given by
\begin{gather*} \Delta(D)= D\otimes D, \ \ \ \Delta(K_i) = K_i\otimes K_i,\\ 
\Delta(X_i^+) = X_i^+\otimes K_i + 1\otimes X_i^+, \ \ \ \Delta(X_i^-) =  X_i^-\otimes 1 + K_i^{-1}\otimes X_i^-,\\
S(K_i) =K_i^{-1}, \ \ \ S(D)= D^{-1}, \ \ \ S(X_i^+) = -X_i^+K_i^{-1}, \ \ \ S(X_i^-) = -K_iX_i^{-},\\
\epsilon(K_i)=1 =\epsilon(D), \ \ \ \epsilon(X_i^\pm)=0.
\end{gather*}

Set $$C=K_0K_\theta,\ \ K_\theta=\prod_{i\in I}K_i^{c_i},$$
where $c_i$ are such that $h_\theta = \sum_{i\in I}c_ih_i$. The quantized enveloping algebra $\bu_q(\lie g)$ is the Hopf subalgebra generated by the elements $X_i^\pm , K_i^{\pm 1}$, $i\in I$.
\subsubsection{The Drinfeld presentation of $\bu_q(\widehat{\lie g})$} An alternate  presentation of the quantum affine algebra was given by Drinfeld. 

The algebra $\bu_q(\widehat{\lie g})$ is isomorphic to the $\mathbb C(q)$-associative algebra  with unit given by generators $c^{\pm 1/2}$, $x_{i,r}^{\pm}$, $k_i^{\pm 1}$, $d^{\pm 1}$ $h_{i,s}$, for $i\in I$, $r,s\in \mathbb Z$ with $s\ne 0$ %and central elements $c^{\pm 1/2}$, 
subject to the following relations:
\begin{gather*}
c^{1/2}c^{-1/2} = 1 = dd^{-1} = k_ik_i^{-1} =  k_i^{-1}k_i, \ \ \ c^{\pm 1/2}\ {\rm are \ central},\\
k_ik_j =k_jk_i, \ \ \  k_ih_{j,r}= \ h_{j,r}k_i, \ \ \ dk_i = k_id, \ \  dh_{i,r}d^{-1} = q^rh_{i,r}\\
k_ix_{j,r}^\pm k_i^{-1} =  q_i^{{}\pm a_{i,j}}x_{j,r}^{{}\pm{}}, \ \ \ \ dx_{i,r}^\pm d^{-1} = q^rx_{i,r}^\pm,\\
[h_{i,r}, h_{j,s}] =  \delta_{r,-s}\frac{1}{r}[r a_{i,j}]_{q_i}\frac{c^r- c^{-r}}{q_j - q_j^{-1}},   \  \ \ [h_{i,r}, x_{j,s}^{{}\pm{}}] =   \pm\frac1r[r a_{i,j}]_{q_i}c^{\mp|r|/2}x_{j,r+s}^{{}\pm{}},\\
x_{i,r+1}^{{}\pm{}}x_{j,s}^{{}\pm{}} -q_i^{{}\pm a_{i,j}}x_{j,s}^{{}\pm{}}x_{i,r+1}^{{}\pm{}}  =  \ q_i^{{}\pm     a_{i,j}}x_{i,r}^{{}\pm{}}x_{j,s+1}^{{}\pm{}} -x_{j,s+1}^{{}\pm{}}x_{i,r}^{{}\pm{}},\\
%[h_{i,n}, h_{j,m}] = \delta_{n,-m} \frac{1}{n}[nc_{ij}]_{q_i}\dfrac{c^n - c^{-n}}{q - q^{-1}}, \ \ & 
 %[x_{i,r}^+ , x_{j,s}^-]= \delta_{i,j} \frac{c^{(r-s)/2}\phi_{i,r+s}^+ - c^{-(r-s)/2}\phi_{i,r+s}^-}{q_i - q_i^{-1}},\\
[x_{i,r}^+ , x_{j,s}^-]= \delta_{i,j} \ \frac{c^{(r-s)/2}\phi_{i,r+s}^+ - c^{-(r-s)/2} \phi_{i,r+s}^-}{q_i - q_i^{-1}},\\
\sum_{\sigma\in S_m}\sum_{k=0}^m(-1)^k \qbinom{m}{k}_{q_i}x_{i,
n_{\sigma(1)}}^{\pm}\ldots x_{i,n_{\sigma(k)}}^{{}\pm{}} 
x_{j,s}^{{}\pm{}} x_{i, n_{\sigma(k+1)}}^{{}\pm{}}\ldots
x_{i,n_{\sigma(m)}}^{{}\pm{}} =0,\ \ \text{if $i\ne j$},
\end{gather*}
for all sequences of integers $n_1,\ldots, n_m$, where $m
=1-a_{i,j}$, $i,j\in I$, $S_m$ is the symmetric group on $m$ letters, and the
$\phi_{i,r}^{{}\pm{}}$ are determined by equating powers of $u$ in
the formal power series
\begin{equation*}
\phi_i^\pm(u) = \sum_{r=0}^{\infty}\phi_{i,\pm r}^{\pm}u^{\pm r} =
k_i^{\pm 1} \exp\left(\pm(q_i-q_i^{-1})\sum_{r=1}^{\infty}h_{i,\pm r}
u^{\pm r}\right).
\end{equation*}
Note that $\phi_{i,\mp r}^{\pm}=0$ for $r>0$.
%The subalgebra $\widehat{\mathcal A_q '}$ of $\widehat{\mathcal A_q}$ generated by $x_{i,r}^\pm ,\  h_{i,s},\  k_{i}^{\pm 1},\  c^{\pm 1/2}$, $i\in I$, $r,s\in \mathbb Z$, $s\neq 0$ is isomorphic to $\bu_q(\widehat{\lie g}')$.
The quantized loop algebra $\bu_q(L(\lie g))$ is the algebra with generators $x_{i,r}^\pm$, $k_i^{\pm 1}$, $h_{i,s}$, $r,s\in\mathbb{Z}$ with $s\neq 0$ and $i\in I$ and the same relations as above where we replace $c^{1/2} $ by 1. Moreover we have a canonical inclusion $\bu_q(\lie g)\hookrightarrow \bu_q(L(\lie g))$ given by mapping $X_i^\pm \to x_{i,0}^\pm$, $k^{\pm}_i\to K^{\pm}_i$ $i\in I$.\\\\
%It will be convenient to consider the elements $\Lambda_{i,r}, i\in I, r\in \mathbb Z$ defined b$$\sum_{r=0}^\infty \Lambda_{i,\pm r}u^r = \exp\left(-\sum_{s=1}^\infty \dfrac{h_{i,\pm s}}{[s]_{q_i}}u^s\right).$$
Explicit formulae for the Hopf algebra structure in terms of these generators are not known. However the following partial information on the coproduct is often enough for our purposes \cite{CP96, Da98}. 

\begin{prop}Let 
$$X^\pm = \sum_{i\in I,\ r\in \mathbb Z}\mathbb C(q)x_{i,r}^\pm, \ \ \ X^\pm(i) = \sum_{j\in I\setminus\{i\},\ r\in \mathbb Z}\mathbb C(q)x_{j,r}^\pm, \ \ \ i\in I.$$
Then 
\begin{enumerit}
\item Modulo $\bu_q(\widehat{\lie g})X^-\otimes \bu_q(\widehat{\lie g})(X^+)^2 + \bu_q(\widehat{\lie g})X^-\otimes \bu_q(\widehat{\lie g})X^+(i)$, we have 
$$\Delta(x_{i,k}^+) = x_{i,k}^+\otimes 1 + k_i\otimes x_{i,k}^+ + \sum_{j=1}^k\phi_{i,j}^+ \otimes x_{i,k-j}^+, \ \ \ k\geq 0,$$
$$\Delta(x_{i,-k}^+) = x_{i,-k}^+\otimes 1 + k_i^{-1}\otimes x_{i,-k}^+ + \sum_{j=1}^{k-1}\phi_{i,-j}^- \otimes x_{i,-k+j}^+, \ \ \ k> 0.$$
\item Modulo $\bu_q(\widehat{\lie g})(X^-)^2\otimes \bu_q(\widehat{\lie g})X^+ + \bu_q(\widehat{\lie g})X^-\otimes \bu_q X^+(i)$, we have
$$\Delta(x_{i,k}^-) = x_{i,k}^-\otimes k_i + 1\otimes x_{i,k}^- + \sum_{j=1}^{k-1} x_{i,k-j}^-\otimes \phi_{i,j}^+, \ \ \ k> 0,$$
$$\Delta(x_{i,-k}^-) = x_{i,-k}^-\otimes k_i^{-1} + 1\otimes x_{i,-k}^- + \sum_{j=1}^{k} x_{i,-k+j}^-\otimes \phi_{i,-j}^-, \ \ \ k\geq 0.$$
\item Modulo  $\bu_q(\widehat{\lie g})X^-\otimes \bu_q(\widehat{\lie g})X^+$, we have
$$\Delta(h_{i,k})= h_{i,k}\otimes 1 + 1 \otimes h_{i,k}, \ \ k\in\mathbb Z\setminus\{0\}.$$
\end{enumerit}\qed 
\end{prop}

\subsection{Representations of quantum algebras}\label{repq}
The classification of finite-dimensional and integrable representations of the quantum algebras is essentially the same as that of the corresponding Lie algebras, once we impose certain restrictions. Thus we shall only be interested in {\em type 1} modules for these algebras; namely we require that the elements $K_i^{\pm 1}$ act semi-simply on the module with eigenvalues in $q^{\mathbb Z}$. The character of such a representation $V$ is given by  $\ch(V) = \sum_{\mu} \dim (V_\mu) e_\mu$ where $V_\mu=\{v\in V: K_i^{\pm 1}v=q_i^{\pm\mu(h_i)}v, \ i\in I\}.$ 
\subsubsection{The irreducible modules} \label{irreducible} Let $\mathcal P_q^+$ (resp. $\mathcal P_q$) be the free abelian monoid  (resp. group) generated by elements $\bomega_{i,a}$ with $a\in\mathbb C(q)^{\times}$.  Clearly we can regard $\mathcal P^+$ as a submonoid of $\mathcal P_q^+$. Define $\weight :\ \mathcal P_q^+\to P^+$ by extending the assignment $\weight\ \bomega_{i,a}=\omega_i$ to a morphism of monoids.

Part (i)   of the next result was proved in \cite{CP91}  while part (ii) was proved in \cite{L2010}, \cite{Ro88}. 
 \begin{thm}\hfill
\begin{enumerate}
\item[(i)] There is a bijection $\bomega\to[V_q(\bomega)]$ between elements of $\mathcal P_q^+$ and isomorphism classes of finite-dimensional irreducible representations of $\bu_q(L(\lie g))$.
   Moreover if $\bomega,\bomega'\in\mathcal P_q^+$ then $V_q(\bomega\bomega')$ is a subquotient of $V_q(\bomega)\otimes V_q(\bomega')$.
   \item[(ii)] Given $\lambda\in P^+$ there exists a unique (up to isomorphism) finite-dimensional irreducible $\bu_q(\lie g)$-module $V_q(\lambda)$. Moreover ${\rm ch} V_q(\lambda)={\rm ch}V(\lambda)$
   and for $\lambda,\mu\in P^+$ we have that $V_q(\lambda+\mu)$ is a summand of   $V_q(\lambda)\otimes V_q(\mu)$. An analogous statement holds for positive level integrable representations of $\bu_q(\widehat{\lie g})$. 
   
   \hfill\qedsymbol
\end{enumerate}
\end{thm}
\begin{rem} \begin{enumerate}
    \item It is worth emphasizing that  in general ${\rm{ch}} V_q(\bomega)\ne{\rm{ch}} V(\bomega)$  for $\bomega\in\mathcal P^+$.
    \item The elements $\bomega_{i,a}\in\cal P_q^+$ are called the fundamental weights and the associated representations are called fundamental representations.
\end{enumerate} 
\end{rem}
\subsubsection{The category $\cal F_q$}\label{prime}
 We shall be interested primarily in the category $\cal F_q$ of finite-dimensional representations of $\bu_q(L(\lie g))$. The Hopf algebra structure of $\bu_q(L(\lie g))$ ensures that $\cal F_q$ is a monoidal tensor category.  Roughly speaking this means that $\cal F_q$ is an abelian  category which is closed under taking tensor products and duals. However, since the Hopf algebra is not co-commutative it is not true in general that the modules $V\otimes W$ and $W\otimes V$ are isomorphic. One has also to be careful to distinguish between left and right duals say $V^*$ and $^*V$; in one case we have an inclusion of $\mathbb C\hookrightarrow  V\otimes V^*$ and in the other case a projection $^*V\otimes V\to\mathbb C\to 0$. The tensor product defines a ring structure on the Grothendieck group of this category. A very interesting fact  proved in  \cite{FR99} is that the Grothendieck ring is always commutative.\\\\
 We shall say that an irreducible representation in $\cal F_q$ is prime if it cannot be written in a nontrivial way as a tensor product of two objects of $\cal F_q$. It is trivially true that any irreducible representation is isomorphic to a tensor product of prime representations. It follows that to understand irreducible representations it is enough to understand the prime ones. However, outside $\lie g=\lie{sl}_2$ the classification of prime objects seems to be a very hard and perhaps wild problem. We shall nevertheless, give various examples of families of prime representations in this section and the next, including some coming from the connection with cluster algebras.\\\\
 Our next definition is entirely  motivated by the connection with cluster algebras. Namely we shall say that an object $V$ of $\cal F_q$ is real if $V^{\otimes r}$ is irreducible for all $r\ge 1$. As a consequence of the main result of \cite{He09} it is enough to require $V^{\otimes 2}$ to be irreducible. Again, it is hard to characterize real representations. A well-known example of Leclerc shows in \cite{Lec02} that there are prime representations which are not real.\\\\
 The notion of prime and real objects can obviously be defined for the finite-dimensional module category of any Hopf algebra, and in particular for irreducible  finite-dimensional representations of $\lie g$ and  $L(\lie g)$. For $\lie g$ it is an exercise to prove that the representations $V(\lambda)$, $\lambda\in P^+$ are prime and not real if $\lambda\ne 0$.
 In the case of $L(\lie g)$ it then follows from Proposition \ref{levelzero} that the prime irreducible representations are precisely $\ev_a V(\lambda)$, $\lambda\in P^+$, $a\in\mathbb C$. Moreover, it also follows that these representations are never real if $\lambda\neq 0$.  So these notions are uninteresting in these examples.\\\\
 As in the case of $L(\lie g)$ the objects of $\cal F_q$ are not completely reducible. In fact these categories behave more like the category $\cal O$ for simple Lie algebras and some of these similarities are explored in this article.
 \subsubsection{Representations of quantum loop $\lie{sl}_2$}\label{sl2} In this case the study of the irreducible objects in $\cal F_q$ is well-understood and we briefly review the main results from \cite{CP91}.  Given $r\in\mathbb Z_+$ and $a\in\mathbb C(q)^\times$ let $$\bomega_{1,a,r}=\bomega_{1,aq^{r-1}}\bomega_{1,aq^{r-3}}\cdots\bomega_{1, aq^{-r+1}}.$$ 
 Note that $\bomega_{1,a,0}$ is the unit element of the monoid $\mathcal{P}_q^+$ for all $a\in \mathbb C(q)^\times$.
 Then $$ V_q(\bomega_{1,a,r})\cong_{\bu_q(\lie{sl}_2)}V_q(r\omega_1).$$ 
 Moreover  $$V_q(\bomega_{1,a,r})\otimes V_q(\bomega_{1,b,s})\cong V_q(\bomega_{1,a,r}\bomega_{1,b,s})\iff ab^{-1}\notin \{q^{\pm(r+s-2p)}: 0\le p<\min\{r,s\}\}.$$  In particular, the modules $V_q(\bomega_{1,a,r})$ are prime and real.
If $ab^{-1}= q^{\pm(r+s-2p)}$ for some $0\le p<\min\{r,s\}$ then we have a non-split short exact sequence,
$$0\to V_1\to V_q(\bomega_{1,a,r})\otimes V_q(\bomega_{1,b,s})\to V_2\to 0,$$
where 
$$ V_1\cong V_q(\bomega_{1,aq^{(r-p)},p})\otimes V_q(\bomega_{1,bq^{(p-r)},r+s-p}),\ \  \ V_2\cong V_q(\bomega_{1,aq^{-p-1},r-p-1})\otimes V_q(\bomega_{1,aq^{p+1},s-p-1})$$ if $ ab^{-1}=q^{-(r+s-2p)}$  while if $ab^{-1}=q^{(r+s-2p)}$  then
$$V_1\cong V_q(\bomega_{1,aq^{p+1},r-p-1})\otimes V_q(\bomega_{1,aq^{-p-1},s-p-1}), \ \  V_2\cong V_q(\bomega_{1,aq^{(p-r)},p})\otimes V_q(\bomega_{1,bq^{(r-p)},r+s-p}).$$
Any irreducible module in $\cal F_q$ is isomorphic to a tensor product of representations of the form $V_q(\bomega_{1,a,r})$, $r\in\mathbb Z_+$, $a\in\bc(q)^{\times}$. More precisely, if $\bomega\in\cal P_q^+$ then $$V(\bomega)\cong V_q(\bomega_{1,a_1,r_1})\otimes \cdots\otimes V_q(\bomega_{1,a_k,r_k}),$$ for a unique choice of $k\ge 1$ and pairs $(a_s,r_s)$, $1\le s\le k$ satisfying $a_sa_m^{-1}\ne q^{\pm(r_s+r_m-2p)}$ for any $0\le p<\min\{r_s,r_m\}$, $1\le s\ne m\le k$.
\subsubsection{Local Weyl modules for quantum loop algebras}\label{locweyl}
  We identify the monoid $\mathcal P_q^+$ with the monoid consisting of $I$-tuples of polynomials $(\pi_i(u))_{i\in I}$, $\pi_i(u)\in \mathbb C(q)[u]$ $\pi_i(0)=1$, via 
  $$\bomega_{i,a}\mapsto (1-\delta_{i,j}au)_{j\in I}.$$
 
For $\bomega\in\mathcal P_q^+$ the local Weyl module $W_q(\bomega)$ is  the $\bu_q(L(\lie g))$-module generated by an element $v_\bomegas$ with relations
 $$x_{i,r}^+ v_{\bomegas} = 0 = (x_{i,0}^-)^{\deg\pi_i(u)+1}v_\bomegas, \ \ \   \phi_{i,r}^\pm v_{\bomegas} = \gamma_{i,r}^\pm v_\bomegas,\ \ r\in\mathbb{Z}$$
 where $\gamma_{i,r}^\pm \in \mathbb C(q)$ are defined by  
 $$ \sum_{r=0}^\infty \gamma_{i,\pm r}^{\pm}u^{\pm r} = q_i^{\deg \pi_i}\frac{\pi_i(q_i^{-1} u)}{\pi_i(q_iu)},
 \ \ \ \ \bomega = (\pi_i(u))_{i\in I}.$$

 The following was proved in \cite{CP01}.
 
 \begin{prop} Let $\bomega\in\mathcal P_q^+$. Then   $\dim W_q(\bomega)<\infty$ and $W_q(\bomega)$ has a unique irreducible quotient $V_q(\bomega)$. If $\bomega'\in\mathcal P_q^+$ then the  modules $W_q(\bomega)$ and $W_q(\bomega')$ are isomorphic if and only if  $\bomega=\bomega'$. The module $V_q(\bomega\bomega')$ is a subquotient of $W_q(\bomega)\otimes W_q(\bomega')$. \hfill\qedsymbol
 \end{prop}
 \subsubsection{The fundamental local Weyl modules}\label{fundloc}  It was proved in \cite{C01}  that $W_q(\bomega_{i,a})\cong V_q(\bomega_{i,a})$. In general it is not true that local Weyl modules are irreducible but we will discuss conditions for these later in the paper. 
 \subsubsection{Local Weyl modules and tensor products}\label{weyltp}  Suppose that $\bomega,\bomega'\in\mathcal P_q^+$ and let $M$, $M'$ be any quotient of $W_q(\bomega)$ and $W_q(\bomega')$ respectively. Let   $v_\bomegas$ and $v_{\bomegas'}$ also denote the images of these elements in $M$ and $M'$. Then using the formulae for comultiplication one can prove  that we have the following sequence of surjective maps:$$W_q(\bomega\bomega')\twoheadrightarrow \bu_q(L(\lie g))(v_{\bomegas}\otimes v_{\bomegas'})\twoheadrightarrow V_q(\bomega\bomega').$$
 In particular  $V_q(\bomega\bomega')$ is a subquotient of $M\otimes M'$.

Suppose that $\bomega=\bomega_{i_1,a_1}\cdots\bomega_{i_k,a_k}\in\cal P_q^+$; the discussion so far establishes the existence 
of a map of $\bu_q(L(\lie g))$-modules $$\phi_{\bomega}: W_q(\bomega)\to V_q(\bomega_{i_1,a_1})\otimes\cdots\otimes V_q(\bomega_{i_k,a_k}).$$
If we assume that $a_1,\dots ,a_k$ are such that 
$a_j/a_s\notin q^{\mathbb N}$ for all $1\le j<s\le k$, then 
the results of  \cite{AK97,Ch01, VV02} show that this map is surjective.
It was conjectured in \cite{CP01} that $\phi_\bomega$  is an isomorphism. Clearly to prove the conjecture it suffices to establish an equality of dimensions. This equality was proved in that paper for $\lie{sl}_2$. In the general case the conjecture was established through the work of \cite{CL06, FoL07, Na11}. We will discuss this further in the next section. \\\\
Assuming  from now on that $\phi_\bomega$ is an isomorphism we study the question of the irreducibility of $W_q(\bomega)$.  
A sufficient condition for $W_q(\bomega)$ to be irreducible is to require that $a_j/a_s\notin q^{\mathbb Z}$ for all $1\le j\ne s\le k$; this was known through the work of \cite{AK97,Ch01}.   A  precise statement was given in \cite[Corollary 5.1]{Ch01} when $\lie g$ is of classical type and for some of the exceptional nodes. The following result summarizes the results of \cite{Ch01} for classical cases.
\begin{thm} 
\begin{enumerate}
\item[(i)] Suppose that $\lie g$ is of classical type and $i,j\in I $ and $a,b\in\mathbb C(q)^{\times}$. Then $$ab^{-1}\notin q^{\pm S(i,j)}\implies V_q(\bomega_{i,a}\bomega_{j,b})\cong V_q(\bomega_{i,a})\otimes V_q(\bomega_{j,b})\cong W_q(\bomega_{i,a}\bomega_{j,b}),$$ where $S(i,j)$ is given as follows: \\

 \noindent If $\lie g$ is of type $A_n$: $$S(i,j)=\{2+2k-i-j: \max\{i,j\}\le k\le \min\{i+j-1, n\}\}.$$
 
 \vspace{0,2cm}
 
\noindent If $\lie g$ is of type $B_n$, and $\alpha_1$ is short: \vspace{0,15cm}

    \begin{itemize}
    \item $S(1,1) = \{4k-2 : 1\leq k\leq n\}$, \vspace{0,15cm}
    \item $S(i,1) = S(i,1)= \{4k - 2i + 1 : \ i\leq k\leq n\}$, \vspace{0,15cm}
    \item $S(i,j) = \{4+4k-2i-2j : \ \max\{i,j\}\leq k\leq n\}\cup \{4k-2-2|j-i|: \ \max\{i,j\}\leq k\leq n\}$, $i,j>1$.
\end{itemize}

 \vspace{0,2cm}
 
\noindent If $\lie g$ is of type $C_n$, and $\alpha_1$ is long: \vspace{0,15cm}
\begin{itemize}
    \item $S(1,1) = \{2k+2 : \ 1\leq k\leq n\}$\vspace{0,15cm}
    \item $S(i,1) = S(1,i)= \{2k-i+3 : 1\leq k\leq n\}$, $i>1$,\vspace{0,15cm}
    \item $S(i,j) = \{2+2k - i- j : \ \max\{i,j\}\leq k\leq n\}\cup \{2+2k-|i-j|: \ \max\{i,j\}\leq k\leq n\}$, $i,j>1$.
\end{itemize}

\vspace{0,2cm}

\noindent If $\lie g$ is of type $D_n$, and $1$ and $2$ denote the spin nodes:
\vspace{0,15cm}

\begin{itemize}
    \item $S(1,1) = S(2,2) = \{2k-2: 2\leq k\leq  n, \ \ k\equiv 0 \mod 2\}$\vspace{0,15cm}
    \item $S(1,2)=S(2,1) = \{2k-2 : 3\leq k\leq n, \ \ k\equiv 1 \mod 2\}$\vspace{0,15cm}
    \item $S(1,j) = S(2,j)= \{2k-j : j\leq k\leq n\} = S(j,2)= S(j,1)$, \ \ $j\geq 3$.\vspace{0,15cm}
    \item $S(i,j) = \{2+2k-i-j : \ \max\{i,j\}\leq k\leq n\}\cup \{-2+2k-|i-j|: \ \max\{i,j\}\leq k\leq n\}$, $i,j\geq 3$. 
\end{itemize}
\vspace{0,2cm}
\item[(ii)]  Given $\bomega=
\bomega_{i_1,a_1}\cdots\bomega_{i_k,a_k}
$ the module $W_q(\bomega)$ is irreducible if $$a_ra_s^{-1}\notin q^{\pm S(i_r,i_s)},\ \ 1\le r<s\le k.$$
\end{enumerate}
\hfill\qedsymbol
\end{thm}

\begin{rem} More recently an alternative approach to describing the set $S(i,j)$ was given in \cite[Theorem 2.10]{Fu22} and \cite[Section 6]{fuoh21} by relating it to the poles of the universal $R$-matrix. It is nontrivial to see that those conditions are equivalent to the explicit description given in the preceding theorem.
\end{rem}

\subsubsection{$\bold A$-forms and classical limits}\label{aform}
Let $\bold A=\mathbb Z[q,q^{-1}]$ and define $\bu_\ba(\widehat{\lie g})$ to be the $\bold A$-subalgebra of $\bu_q(\widehat{\lie g})$ generated by the elements $(X_i^\pm)^r/[r]_q!$, $i\in\widehat I$. Then $$\bu_q(\widehat{\lie g})\cong \bu_{\bold  A}(\widehat{\lie g})\otimes_{\bold A} \mathbb C(q).$$ For $\epsilon\in\bc^\times$ we let $\bc_\epsilon$ be the $\bold A$-module obtained by letting $q$ act as $\epsilon$ and set $$\bu_\epsilon(\widehat{\lie g})= \bu_{\bold  A}(\widehat{\lie g})\otimes_{\bold A} \mathbb C_\epsilon.$$ The algebra $\bu(\widehat{\lie g})$ is isomorphic to the quotient of $\bu_1(\widehat{\lie g})$ by the ideal generated by $K_i-1, D-1,\ i\in \widehat I$. Similar assertions hold for $\bu_q(\lie g)$  and $\bu_q(L(\lie g))$ as well. 
Part (i) of the following was proved in \cite{L2010}, \cite{Ro88} while part (ii) was proved in \cite{CP01}.
 \begin{thm}\hfill
 \begin{enumerate}
\item [(i)] Suppose that $\lambda\in\widehat P^+$ and $\epsilon\in\mathbb C^\times $. There exists a  $\bu_\bold A(\widehat{\lie g})$-submodule $V_\ba(\lambda)$ of $V_q(\lambda)$ such that $$V_q(\lambda)\cong V_\ba(\lambda)\otimes_\ba\bc(q).$$
In particular $V_\epsilon(\lambda)= V_\ba(\lambda)\otimes_\ba\bc_\epsilon$ is a module for $\bu_\epsilon(\widehat{\lie g})$, and if $\epsilon=1$ or if $\epsilon$ is not a root of unity then we have ${\rm{ch}} V_\epsilon(\lambda)={\rm{ch}} V_q(\lambda).$ Analogous statements hold for the $\bu_q(\lie g)$ representations $V_q(\lambda)$,  $\lambda\in P^+$.\\
\item[(ii)] Let $\mathcal P_\ba^+$ be the submonoid of $\mathcal P^+$ generated by elements $\bomega_{i,a}$, $a\in\bold A$ and let $\bomega\in \mathcal P_\ba^+$. Then $W_q(\bomega)$
admits a $\bu_\ba(L(\lie g))$-submodule $W_\ba(\bomega)$ and $$W_q(\bomega)\cong W_\ba(\bomega)\otimes_\ba\bc(q).$$ In particular $W_\epsilon(\bomega)=W_\ba(\bomega)\otimes _\ba\bc_\epsilon$ is a $\bu_\epsilon(L(\lie g))$-module, and ${\rm{ch}} W_\epsilon(\bomega)= W_q(\bomega)$ if $\epsilon=1$ or if $\epsilon$ is not a root of unity. 
If $M_q$ is  any $\bu_q(L(\lie g))$-module quotient of $W_q(\bomega)$ let  $M_{\bold A}$ be the image of $W_{\bold A}(\bomega)$. Then,
$$M_q\cong M_{\bold A}\otimes_{\bold A} \mathbb C(q),\ \ M_\epsilon\cong M_{\bold A}\otimes_{\bold A} \mathbb C_\epsilon,$$ and $M_\epsilon$ is a canonical quotient of $W_\epsilon(\bomega)$.
\end{enumerate}
\hfill\qed
 \end{thm}
 \begin{rem} \mbox{}
 \begin{itemize}
\item The modules  $V_1(\lambda)$, $\lambda\in \widehat P^+$ and $V_1(\bomega)$, $\bomega\in \mathcal P^+_\ba$ are modules for the universal enveloping algebra of $\bu(L(\lie g))$ and are called the classical limit of $V_q(\lambda)$ and $V_q(\bomega)$ respectively.
\item It is worth noting again,  that in part (i) of the theorem the module $V_\epsilon(\lambda)$ is irreducible for $\bu_\epsilon(\widehat{\lie g})$ if $\epsilon$ is not a root of unity. This is false in part (ii).
 \end{itemize}
\end{rem}

\section{Classical and Graded limits}\label{classgrad}
In this section we discuss various well-studied families of finite-dimensional representations of quantum affine algebras. We review the literature on the  presentation of these modules, their classical limits and the closely related  {\em  graded limits}. 
\subsection{Classical and Graded Limits of the Quantum Local Weyl modules}\label{locweyldim} 
\subsubsection{Relations in the classical limit} \label{classicallt} Suppose that $\bomega\in\mathcal P^+_{\bold A}$ and let $V_q(\bomega)$ be the unique irreducible quotient of $W_q(\bomega)$ (see Section \ref{locweyl}).  Since $W_q(\bomega)$ is  finite-dimensional the corresponding classical limits (see Theorem \ref{aform}) $W_1(\bomega)$ and $V_1(\bomega)$ are finite-dimensional modules for $L(\lie g)$. Let ${\bar v_{\bomega}}=v_\bomega\otimes 1\in W_1(\bomega)$. The following was proved in \cite{CP01}.

\begin{lem} Suppose that $\bomega=\bomega_{i_1,a_1}\cdots \bomega_{i_k,a_k}\in\cal P^+_{\bold A}$. The following relations hold in $W_1(\bomega)$: $$h\bar{v}_\bomega=\weight\ \bomega(h)\bar{v}_\bomega,\ \  \left(h\otimes (t-a_1(1))\cdots (t-a_k(1))\bc[t, t^{-1}]\right){\bar v_{\bomega}}=0,$$ $$ (x_\alpha^+\otimes\mathbb C[t,t^{-1}]){\bar v_{\bomega}}=0,\ \  h\in\lie h,\ \ \alpha\in  R^+.$$
\hfill\qedsymbol
\end{lem}
\begin{rem} With a little more work one can actually prove that there exists an integer $N\in\mathbb N$ such that  $W_1(\bomega)$ is a module for the truncation of $L(\lie g)$ by the polynomial $(t-a_1(1))^N\cdots (t-a_k(1))^N$.
\end{rem}
\subsubsection{Graded Limits: the modules $V_{\loc}$} \label{gradedlt} {\em In the rest of this section we shall restrict our attention to the submonoid $\cal P^+_{\mathbb Z}$ which is generated by the elements $\bomega_{i,q^r}$ for $i\in I$ and $r\in\mathbb Z$.}\\\\
 In this case the results of Section \ref{classicallt} imply that $W_1(\bomega)$ is a module for the truncation of $L(\lie g)$ at $(t-1)^N$ for some $N$ sufficiently large.  Using the isomorphism of Lie algebras $$\lie g\otimes\frac{\mathbb C[t,t^{-1}]}{(t-1)^N}\cong \lie g\otimes\frac{\mathbb C[t]}{(t-1)^N}\cong \lie g\otimes\frac{\mathbb C[t]}{(t^N)}$$ we see that we can regard $W_1(\bomega)$ as a module for the truncation of  $\lie g[t]$ at $t^N$. We shall denote this module by $W_{\loc}(\bomega)$. We call this the {\em graded limit} of $W_q(\bomega)$. In fact if $M$ is any quotient of $W_q(\bomega)$ we can define a corresponding module $M_{\loc}$ for $\lie g[t]$ using the isomorphisms of Lie algebras and call this the graded limit of $M$.\\\\
 This terminology of course requires justification. Recall that the action of $d$ on $\lie g[t]$ defines a $\mathbb Z_+$-grading on it: the $r$-th graded piece is $\lie g\otimes t^r$. The adjoint action of $d$   on  $\bu(\lie g[t])$ also gives a $\mathbb Z_+$-grading. Hence one can define the notion of a graded $\lie g[t]$-module $V$ to be one which admits a compatible $\mathbb Z$-grading namely: $$V=\bigoplus_{s\in\mathbb Z}V[s],\ \ (\lie g\otimes t^r)V[s]\subseteq V[r+s].$$ The general belief is that when $M$ is a quotient of $W_q(\bomega)$ with $\bomega\in\cal P^+_{\mathbb Z}$ then $M_{\loc}$ is a graded $\lie g[t]$-module. This is far from clear in general and is hard to prove even in specific cases. In the rest of the section we will discuss certain families of modules where the corresponding graded limit is in fact a graded $\lie g[t]$-module.
 
 \begin{rem} In the discussion that follows we shall see that the classical or graded limit depends only on $\weight~\bomega$. So there is a substantial loss of information when we go to the limits. However, the character and the underlying $\bu_q(\lie g)$-module is the same and this is one reason for our interest in this study.
 
 \end{rem}
 \subsubsection{Kirillov--Reshetikhin modules} \label{krq}
 We begin by discussing this particular family of modules since this was essentially the motivation for the interest in graded limits. 
 
Given $i\in I$, $r\in\bn$ and $s\in\bz$ set  $$\bomega_{i,s,r}=\bomega_{i,q_i^{s+r-1}}\bomega_{i, q_i^{s+r-3}}\cdots\bomega_{i,q_i^{s-r+1}}.$$ 
Notice that in the case $i=1$ these elements were introduced in Section \ref{sl2} in the case of $\lie{sl}_2$ where they were denoted as $\bomega_{1,q^s,r}$ since we were working in a more general situation.
The corresponding irreducible $\bu_q(L(\lie g))$-module  is called a Kirillov--Reshetikhin module. This is because of an important conjecture that they had made; they predicted the existence of certain modules for the quantum loop algebra with a specific decomposition as  $\bu_q(\lie g)$-modules (see also \cite{HKOTY99}).  In \cite{C01} it was proved that the conjectured modules were of the form $V_q(\bomega_{i,s,r})$ for all classical Lie algebras and for some $i\in I$ in the exceptional cases.  Moreover the following  presentation was given (see Corollary 2.1 of \cite{C01}) for the module $V_1(\bomega_{i,s,r})$.
 \begin{thm} The $L(\lie g)$-module $V_1(\bomega_{i,s,r})$ is generated by an element $v_{i,s,r}$ with relations:
$$x_\alpha^+ v_{i,s,r}=0,\ \ (h\otimes t^k) v_{i,s,r}=r\omega_i(h)v_{i,s,r},\ \ ((x_i^-\otimes t^k)-x_i^-\otimes 1)v_{i,s,r}=0, \ \ (x_\alpha^-\otimes 1)^{r\omega_i(h_\alpha)+1}v_{i,s,r}=0,$$ where $\alpha\in R^+$, $h\in\lie h$ and $k\in\mathbb Z$. \hfill\qedsymbol \end{thm}  Here we have used the fact that $\weight~\bomega_{i,s,r}=r\omega_i$. Notice that these relations are independent of $s$. Moreover, $$ ((h\otimes t^k)-h)v_{i,s,r}=0\implies (h\otimes (t-1)^k){v}_{i,s,r}=0,\ \ \ 0\ne k\in\mathbb Z,$$ 
 and similarly for the third relation in the presentation above. 
 It follows that $V_{\loc}(\bomega_{i,s,r})$ is the $\lie g[t]$-quotient of $W_{\loc}(\bomega_{i,s,r})$ by imposing the additional relation: $(x_i^-\otimes t)v_{\bomegas_{i,s,r}}=0.$
 
 Later, in \cite{CM06},   a more systematic self contained study of these modules was developed and the graded $\lie g$-module decomposition of these modules was calculated.
 One can think of this as a graded version of the Kirillov--Reshetikhin character formula. The results of \cite{CM06} led to the definition of graded limits and more generally resulted in the development of the subject of graded (not necessarily finite-dimensional) representations of $\lie g[t]$. We say more about this study in later sections of the paper. 
\subsubsection{Minimal Affinizations}\label{minaff} The notion of minimal affinizations was introduced and further studied in \cite{Ch95b,CP95a,CP95}. Perhaps the simplest place to explain what  this notion means  is in the case of $\lie{sl}_2$. Since we have only one simple root we denote the generators of $\cal P_q^+$ by $\bomega_{1,a}$. If $\bomega=\bomega_{1,a_1}\cdots \bomega_{1,a_k}\in\cal P_q^+$ then it is not hard to see that there exists a $\bu_q(\lie g)$-module $M$ such that $$V_q(\bomega)\cong_{\bu_q(\lie g)} V_q(k\omega)\oplus M.$$ 
Moreover,
it was shown in  \cite{CP91} that $M\ne 0$ unless $V_q(\bomega)$ is a  Kirillov--Reshetikhin module:
$$M=0\iff \bomega=\bomega_{1,aq^{k-1}}\cdots \bomega_{1,aq^{-k+1}},\ \ a\in\bc(q)^{\times}.$$
Using the results of \cite{CP91} we can also give precise conditions under  which $V_q(\bomega)$ and $V_q(\bomega')$ are isomorphic as $\bu_q(\lie g)$-modules.\\\\
It is natural to ask what  analogs of these results hold in the higher rank case. It was known essentially from the beginning (see \cite{Dri88})  that if $\lie g$ is not of type $A$, there does not exist a corresponding  $\bu_q(L(\lie g))$-module structure on $V_q(\lambda)$. On the other hand it is also clear that there were many pairs $\bomega, \bomega'\in\cal P_q^+$ with $V_q(\bomega)\cong_{\bu_q(\lie g)} V_q(\bomega')$. So this motivated the question: given $\lambda\in P^+$,  is there a \lq\lq smallest\rq\rq~ $\bu_q(\lie g)$-module containing a copy of  $V_q(\lambda)$ which admits an action of the quantum loop algebra. This question can be more formally stated as follows.\\\\
Given $\bomega,\bomega'\in\cal P_q^+$ we say that $V_q(\bomega)$ is equivalent to $V_q(\bomega')$ if they are isomorphic as $\bu_q(\lie g)$-modules. Denote the equivalence class corresponding to $\bomega$ by $[V_q(\bomega)]_{\lie g}$. In particular,  $$[V_q(\bomega)]_{\lie g}=[V_q(\bomega')]_{\lie g}\implies \weight~\bomega=\weight~\bomega'.$$
The converse statement is definitely false, this  is already  the case in  $\lie{sl}_2$.\\\\
Define a partial order on the set of equivalence classes by: $[V_q(\bomega)]_{\lie g}\leq [V_q(\bomega')]_{\lie g}$ if for all $\mu\in P^+$ either 
$$\dim\Hom_{\bu_q(\lie g)}(V_q(\mu), V_q(\bomega))\le  \dim\Hom_{\bu_q(\lie g)}(V_q(\mu), V_q(\bomega'))$$
or there exists $\nu> \mu$ (i.e. $\nu-\mu\in Q^+\backslash\{0\}$) such that 
$$\dim\Hom_{\bu_q(\lie g)}(V_q(\nu), V_q(\bomega))<  \dim\Hom_{\bu_q(\lie g)}(V_q(\nu), V_q(\bomega')).$$
It was proved in \cite{Ch95b} that minimal elements exist in this order and an irreducible representation corresponding to a minimal element was called a minimal affinization.  When $\lie g$ is not of type $D$ or $E$ the explicit  expression for the  elements $\bomega\in\cal P_q^+$ which give a minimal affinization are given by 
\begin{gather*}\bomega=\bomega_{i_1,s_1,r_1}\cdots\bomega_{i_k,s_k,r_k},\ \ i_1<i_2<\cdots <i_k,\\  s_{p+1}-s_p = \epsilon\left( d_{i_p}r_p + d_{i_{p+1}}r_{p+1} + \sum_{j=i_{p}}^{i_{p+1}-1}(d_j - 1 - a_{j,j+1})\right),\ \ 1\le p\le k-1,
\end{gather*}
where   either $\epsilon=1$ for all $p$  or $\epsilon=-1$ for all $p$
\iffalse \begin{gather*}\bomega=\bomega_{i_1,s_1,r_1}\cdots\bomega_{i_k,s_k,r_k},\ \ i_1<i_2<\cdots <i_k,\\  s_{p+1}-s_p=-(c_{i_p}(\lambda)+c_{i_p+1}(\lambda)+\cdots c_{i_{p+1}-1}(\lambda)+c_{i_{p+1}}(\lambda)),\ \ 1\le p\le k-1\\ \ {\rm{or}}\ \  s_{p+1}-s_p= 2(d_{i_p}-d_{i_{p+1}})+(c_{i_p}(\lambda)+c_{i_p+1}(\lambda)+\cdots c_{i_{p+1}-1}(\lambda)+c_{i_{p+1}}(\lambda)),\ \ 1\le p\le k-1\end{gather*}
Here $\bomega_{i,s,r}$ is the element of $\cal P^+$ which was introduced in Section \ref{krq} and $$c_i(\lambda)=\begin{cases}d_i(\lambda(h_i)+\lambda(h_{i+1})+1),\ \ (\alpha_i,\alpha_i)=(\alpha_{i+1},\alpha_{i+1}),\\
d_i\lambda(h_i)+d_{i+1}\lambda(h_{i+1})+2d_{i+1}-1,\ \ (\alpha_i,\alpha_i)\ne (\alpha_{i+1},\alpha_{i+1}).\end{cases} $$\fi
and $\bomega_{i,s,r}$ is the element of $\cal P_q^+$ which was introduced in Section \ref{krq}.
In types $D$ and $E$ the preceding formulae still correspond to minimal affinizations under suitable restrictions. Unfortunately these are far from being all of them; the difficulty lies in the existence of the trivalent node (see \cite{CP96ab,CP96}). The problem of  classifying all the  minimal elements was studied in \cite{Per14} but the full details  are still to appear. \\\\
 The equivalence classes of Kirillov--Reshetikhin modules are clearly minimal affinizations. The following result was conjectured in \cite{M10a} and proved in 
 \cite{NL16a,N13a, Na14a}. It again justifies the use of the term graded limit. We do not state the result in full generality in type $D$ but restrict our attention to the minimal affinizations discussed here. 
 \begin{thm} Assume that $\bomega\in\cal P_q^+$ is as in the discussion just preceding the theorem. Then $V_{\loc}(\bomega)$ is the $\lie g[t]$-module generated by an element $v_\bomega$ with relations:
 $$x_{i}^+v_\bomega=0,\ \ (h_i\otimes t^k)v_\bomega=\delta_{k,0}\weight~\bomega(h_i)v_\bomega\ \ (x_\beta^-\otimes t)v_\bomega= 0\ \ (x_i^-\otimes 1)^{\weight~\bomega(h_i)+1}v_\bomega=0,$$
 for all $i\in I$ and for all $\beta=\sum_{i=1}^ns_i\alpha_i\in R^+$  with $s_i\le 1$.
 \hfill\qedsymbol
  \end{thm}
 The $\lie g$-module decomposition of the minimal affinizations was computed in rank two in \cite{Ch95b}. In the case of $A_n$ the graded limit is irreducible and so its character is just the character of $V_q(\weight~\bomega)$. The $\lie g$-module decomposition in type $B$ and $D_4$ was partially given in  \cite{M10a} and the result in complete generality is in \cite{N13a,Na14a} for types $B,C$ and for certain minimal affinizations in type $D$. Moreover, in \cite{Sam13} Sam proved a conjecture made in \cite{CG11} that the character of minimal affinizations in types $BCD$ are given by a Jacobi--Trudi determinant.

 \subsubsection{Tensor products and Fusion products} Before continuing with our justification for the term graded limit, we discuss the following natural question. Suppose that $M$ and $M'$ are $\bu_q(L(\lie g))$-modules which have a classical (resp. graded) limit. Is it true that $M\otimes M'$ has a classical limit and how does this relate to the tensor product of the classical (resp. graded) limits? The answer to this question is far from straightforward; even if $M\otimes M'$ does have a classical limit it is easy to generate examples 
 where $(M\otimes M')_1$ is not  isomorphic to $M_1\otimes M'_1$ as $L(\lie g)$-modules.
 For instance if we take $\lie g=\lie{sl}_2$ and  $M=V_q(\bomega_{1,q^2})$, $M'=V_q(\bomega_{1,1})$ then the module $M\otimes M'$ is a cyclic $\bu_q(L(\lie{sl}_2))$-module and so the classical limit is a cyclic indecomposable module for $L(\lie{sl}_2)$. However  the tensor product of the  classical limits is  $V(\bomega_1)\otimes V(\bomega_1)$ (equivalently $\ev_1 V(\omega)\otimes \ev_1 V(\omega)$) which is completely reducible.\\\\
 The $\lie g$-module structure however is unchanged in the process. This is because it is known that the process of taking classical limits preserves tensor products for the simple Lie algebras. If we work with the graded limit then again it is false that the graded character of $(M\otimes M')_{\loc}$ is the same as  $M_{\loc}\otimes M'_{\loc}$.
However in many examples the  graded  limit of the tensor product coincides with an operation called the fusion product defined on graded $\lie g[t]$-modules. This notion was introduced in \cite{FL99} and we now recall this construction.
 
Let  $V$ be  a finite-dimensional cyclic $\lie g[t]$-module generated by an element $v$ and for $r\in\bz_+$  set  $$F^rV = \left(\bigoplus_{0\leq s \leq r} \bu(\lie g[t])[s]\right)\cdot v$$ Clearly $F^rV$ is a $\lie g$-submodule of $V$ and we have a finite $\lie g$-module filtration $$0\subseteq F^0V\subseteq F^1V\subseteq\cdots \subseteq F^kV =V,$$  for some $k\in\bz_+$. The associated graded vector space $\gr V$ acquires a graded $\lie g[t]$-module structure in a natural way and is  generated by the image of $v$ in $\gr V$.
Given  a  $\lie g[t]$-module $V $ and $z\in\bc$, let $V^z$ be the $\lie g[t]$-module with action $$(x\otimes t^r) w= (x\otimes (t+z)^r)w,\ \ x\in\lie g,\ \  r\in\bz_+,\ w\in V.$$ Let $V_s$, $1\le s\le p$, be cyclic finite-dimensional $\lie g[t]$-modules with cyclic vectors $v_s$, $1\le s\le p$ and let $z_1,\dots, z_p$ be distinct complex numbers. Then the module $V_1^{z_1}\otimes\cdots\otimes V_p^{z_p}$ is cyclic with cyclic generator $v_1\otimes\cdots\otimes v_p.$ The fusion product $V_1^{z_1}*\cdots*V_p^{z_p}$ is defined to be $\gr V_1^{z_1}\otimes\cdots\otimes V_p^{z_p}.$
For ease of notation we shall use $V_1*\cdots*V_k$ for $V_1^{z_1}*\cdots*V_k^{z_k}$.\\\\
It is conjectured in \cite{FL99}  that under some suitable conditions on $V_s$ and $v_s$, the fusion product is independent of the choice of the complex numbers $z_s$, $1\leq s\leq k$, and this conjecture is verified in many special cases by various people (see for instance \cite{CL06},  \cite{CSVW16},  \cite{FF02}, \cite{FL04}  \cite{FoL07}, \cite{Kedem11}, \cite{Na17}). In all these cases the conjecture is proved by exhibiting a graded presentation of the fusion product which is independent of all parameters. This is much like what we have been doing to justify the use of { \em{graded limit}} and the coincidence is not accidental.  In almost all of these papers the proof of the Feigin--Loktev conjecture involves giving a presentation of the graded limit of certain $\bu_q(L(\lie g))$-modules.

\subsubsection{A presentation of $W_{\loc}(\bomega)$} \label{locweyldef} We return to our discussion in Section \ref{weyltp}. Recall that we had discussed  that given $\bomega\in\cal P_q^+$ we can write $\bomega=\bomega_{i_1,a_1}\cdots\bomega_{i_k, a_k}$ so that there is a surjective map of $\bu_q(L(\lie g))$-modules $$W_q(\bomega)\to V_q(\bomega_{i_1,a_1})\otimes \cdots\otimes V_q(\bomega_{i_k,a_k})\to 0.$$
It has been conjectured in \cite{CP2001} that \begin{equation}\label{equaldim}\dim W_q(\bomega)=\prod_{s=1}^k\dim V_q(\bomega_{i_s, a_s})= \prod_{s=1}^k\dim W_q(\bomega_{i_s, a_s}) ,\end{equation} where the second equality is a consequence of Section \ref{fundloc}.
Since the dimension is unchanged when passing to the graded limit it suffices to prove that  \begin{equation}\label{equaldim}\dim W_{\loc}(\bomega)=\prod_{s=1}^k\dim V_{\loc}(\bomega_{i_s, a_s}).\end{equation}
Using Lemma \ref{classicallt} and the discussion in Section \ref{gradedlt} we see that $W_{\loc}(\bomega)$ is the quotient of the module $\tilde W_{\loc}(\weight\bomega)$ which is  generated as a $\lie g[t]$-module by an element $w_{\bomega}$ with defining relations:
$$(h\otimes t^{}\bc[t])w_\bomega=0,\ \  hv_\bomegas=\weight\bomega(h)w_\bomegas,\ \ x_i^+w_{\bomegas}=0,\
(x_i^-)^{\weight\bomegas(h_i)+1}v_{\bomegas}=0.$$  Notice that by Theorem \ref{krq} we know that $$\tilde W_{\loc}(\omega_i)\cong W_{\loc}(\bomega_{i,a}).$$ Choosing distinct scalars $z_1,\dots,z_k$ consider the fusion product $\tilde W_{\loc}(\omega_{i_1})^{z_1}*\cdots  * \tilde W_{\loc}(\omega_{i_k})^{z_k}$. It is not too hard to prove that this module is a quotient of $\tilde W_{\loc} (\weight\bomega)$. We get $$\dim \tilde W_{\loc}(\weight\bomega)\ge\dim W_{\loc}(\bomega)\ge \prod_{s=1}^k\dim \tilde W_{\loc}(\omega_{i_s}).$$
 The following result was established in \cite{CP01} for $\lie{sl}_2$. Using this the result was established in \cite{CL06} in the case of $\lie{sl}_{n+1}$ where a Gelfand--Tsetlin type basis was also given for $W_{\loc}(\bomega)$. These bases were further studied in \cite{RRV14, RRv18}. In \cite{FoL07} the theorem was proved for simply-laced Lie algebras. Finally in \cite{Na11} the result was established for non-simply laced types.
\begin{thm} We have an isomorphism of $\lie g[t]$-modules:
$$\tilde W_{\loc}(\weight \bomega)\cong W_{\loc}(\bomega)\cong \tilde W_{\loc}(\omega_{i_1})^{z_1}*\cdots  * \tilde W_{\loc}(\omega_{i_k})^{z_k},\ \ \weight\bomega=\omega_{i_1}+\cdots+\omega_{i_k}.$$
\hfill\qedsymbol
\end{thm}
Clearly this theorem establishes  the conjecture in \cite{C01} and  also the conjecture of Feigin--Loktev for this particular family of modules.
\begin{rem}
Although the preceding theorem is uniformly stated the methods of proof are very different. In \cite{CL06,CP01} the proof goes by writing down a basis and then doing a dimension count. In \cite{FoL07} the proof proceeds by showing that  $\tilde W_{\loc}(\weight\bomega)$ is isomorphic to  a stable Demazure module in a level one representation of the affine Lie algebra (see Section \ref{levelr} for the relevant definitions). This isomorphism fails in the non-simply laced case. Instead it is proved in \cite{Na11} that the module has a flag by stable level one Demazure modules and this plays a key role in the proof. We return to these ideas in the later sections of this paper.
\end{rem}

\subsubsection{Tensor products of Kirillov--Reshetikhin modules} \label{krtp} It was proved in \cite{Ch01} that  the tensor products of Kirillov--Reshetikhin  modules $V_q(\bomega_{i_1,s_1,r_1})\otimes \cdots\otimes V_q(\bomega_{i_k,s_k,r_k})$ is irreducible as long as $s_i-s_p$, $1\le i\ne p\le k$ lie outside a finite set. A precise description of this set was also given in that paper when  $\lie g$ is classical. Set $$\bold V=V_q(\bomega_{i_1,s_1,r_1})\otimes \cdots\otimes V_q(\bomega_{i_k,s_k,r_k}),\ \ \lambda=\sum_{s=1}^kr_s\omega_{i_s}.$$
We now discuss the results of \cite{Na17} on the structure of $\bold V_{\loc}$. Thus, let
$\tilde{\bold V}_{\loc}$ be the $\lie g[t]$-module generated by a vector $v$ satisfying the relations: 
\begin{align*}
    \lie n^+[t]v= 0 = \left(\lie h\otimes t\mathbb C[t]\right) v, \ \ \ hv= \lambda(h)v, \ \ h\in \lie h,\\
    \left(F_i(z)^r\right)_s v=0, \ \ i\in I, \ r>0, \  s<-\sum_{p:\ i_p=i}\min\{r,r_p\},
\end{align*}
where  $(F_i(z)^r)_s$ denotes the coefficient of $z^s$ in the $r$-th power of  $$F_i(z) = \sum_{m=0}^\infty (x_i^-\otimes t^m)z^{-m-1}\in \bu(\lie g[t])[[z^{-1}]].$$
The following is the main result of \cite{Na17}.
\begin{thm}
We have an isomorphism of graded $\lie g[t]$-modules $$\bold V_{\loc}\cong V_{\loc}(\bomega_{i_1,s_1,r_1})^{z_1}* \cdots*V_{\loc}(\bomega_{i_k,s_k,r_k})^{z_k}\cong \tilde{\bold V}_{\loc}.$$\hfill\qedsymbol
\end{thm}
 Again, the conjecture of Feigin--Loktev for this family of modules is a consequence of this presentation. The proof of the Feigin--Loktev conjecture  when $r_1=r_2=\cdots=r_k$ and $\lie g$ simply-laced was proved earlier in \cite{FoL07} by identifying the fusion product with a $\lie g$-stable Demazure module.
In general the connection with Demazure modules or the existence of a Demazure flag (as in the case of local Weyl modules) is not known.\

\subsubsection{Monoidal categorification and HL-modules}\label{hlm}
Our final example of graded limits comes from the work of David Hernandez and Bernard Leclerc on monoidal categorification of cluster algebras.  We refer the reader to \cite{Wil14} for a quick introduction to cluster algebras. For the purposes of this article it is enough for us to recall that a cluster algebra is a commutative ring with certain distinguished generators called {\em cluster variables} and  certain algebraically independent subsets of cluster variables called {\em clusters}. Monomials in the cluster variables belonging to a cluster are called {\em cluster monomials}. There is also an operation called mutation; this is a way to produce a new cluster by replacing exactly one element of the original cluster by another cluster variable.

The remarkable insight of Hernandez--Leclerc was to relate these ideas 
to the representation theory of quantum affine algebras associated to simply-laced Lie algebras. 
Broadly speaking they prove that the Grothendieck ring of a  suitable tensor subcategory admits the structure of a cluster algebra.  A cluster variable is a prime real representation  in this category (see Section \ref{prime} for the definitions) and we call these the HL-modules. Suppose that $V,V'$ are irreducible modules in this subcategory. Assume that their isomorphism classes correspond to cluster variables which belong to the same cluster. Then $V\otimes V'$ is an irreducible module. The operation of mutation in this language corresponds to the Jordan--Holder decomposition of the corresponding tensor product. \\\\
We now give one specific example of their work and relate it to our study of graded limits. We assume that $\lie g$ is of type $A_n$. Let $\kappa: \{1,\dots,n\}\rightarrow \mathbb{Z}$ be a height function; namely it satisfies  $|\kappa(i+1)-\kappa(i)|= 1$ for $1\leq i\leq n$. Let $\cal P^+_\kappa$ be the submonoid of $\cal P^+_q$ generated by elements $\bomega_{i, q^{\kappa(i)\pm 1}}$, $i\in I$. Let $\cal F_\kappa$ be the full subcategory of $\cal F_q$ consisting of finite-dimensional $\bu_q(L(\lie g))$-modules whose Jordan--Holder constituents are isomorphic to $V_q(\bomega)$ for some $\bomega\in\cal P_\kappa^+$. It was shown in \cite{HL10} that $\cal F_\kappa$ is closed under taking tensor products and that its Grothendieck ring has the structure of a cluster  algebra of type $A_n$. The following result was proved in \cite{HL10} when $\kappa(i)= i\mod 2$,  in \cite{HL13} when $\kappa(i)=i$ and in complete generality in \cite{BC19a}. 
\begin{thm} Suppose that $V_q(\bomega)$ is a prime real object of $\cal F_\kappa$. Then $\bomega$ must be one of the following:  \begin{gather*}\bomega_{i,q^{\kappa(i)\pm 1}},\ \  \bomega_{i,q^{\kappa(i)+1}}\bomega_{i,q^{\kappa(i)-1}},\ \ i\in I,\\
\bomega_{i,a_1}\bomega_{i_2,a_2}\cdots \bomega_{i_{k-1},a_{k-1}}\bomega_{j,a_k},\ 1\leq i<j\leq n,
\end{gather*}
where $i_2<\cdots<i_{k-1}$ is an ordered enumeration of  $\{p: i<p<j,\  \kappa(p-1)=\kappa(p+1)\}$ and $a_1=q^{\kappa(i)\pm 1}$ if $\kappa(i+1)=\kappa(i)\mp 1$ and $a_s=q^{\kappa(i_s)\pm 1}$ if $\kappa(i_s)=\kappa(i_s-1)\pm 1$ for $s\geq 2$. Conversely the irreducible representation associated to any $\bomega$ as above is a real prime object of $\cal F_\kappa$. \hfill\qedsymbol
\end{thm}
\subsubsection{Graded Limits of HL-modules in $\cal F_\kappa$} Continue to assume that $\lie g$ is of type $A_n$ and for $1\le i\le j\le n$ set $\alpha_{i,j}=\alpha_i+\cdots+\alpha_j\in R^+$. It  follows from the discussion in Section \ref{fundloc} that $V_{\loc}(\bomega_{i,q^{\kappa(i)\pm 1}})\cong W_{\loc}(\bomega_{i,q^{\kappa(i)\pm 1}})$. The discussion in Section \ref{krq} gives a presentation for $V_{\loc}(\bomega_{i,q^{\kappa(i)+1}}\bomega_{i,q^{\kappa(i)-1}})$ since this is a special example of a Kirillov--Reshetikhin module.
The following was proved in \cite{BCMo15} and shows  that the graded limits of HL-modules are indeed graded.
\begin{thm}
Suppose that $\lie g$ is of type $A_n$ and $\bomega=\bomega_{i,a_1}\cdots\bomega_{j,a_k}\in\cal P^+$ is as in Theorem \ref{hlm}. Then $V_{\loc}(\bomega)$ is the quotient of $W_{\loc}(\bomega)$ by the submodule generated by the additional relations:
$$(x_\alpha^-\otimes t)w_{\bomega}=0,\ \ \alpha\in \{\alpha_{i,i_2},\alpha_{i_2,i_3},\cdots,\alpha_{i_{k-1},j}\}.$$\hfill\qedsymbol
\end{thm}
We remark that the result in \cite{BCMo15} is more general in the sense that it gives a presentation of the graded limit of the tensor product of an HL-module with the Kirillov--Reshetikhin modules in this category. Here again the result shows that tensor products specialize to fusion products. A problem that has not been studied so far is to understand the graded limit of a tensor product of $V_q(\bomega)\otimes V_q(\bomega')$ for an arbitrary pair $\bomega,\bomega'\in\cal P^+_\kappa$ and the connection with the fusion product of the graded limits of $V_q(\bomega)$ and $V_q(\bomega')$. \\\\
The graded characters of the limits of HL-modules have been studied in \cite{BCSW21} and \cite{BK20b} in different ways. In the first paper a character formula was given as an explicit linear combination of Macdonald polynomials. In \cite{BK20b} the authors studied the $\lie g$-module decomposition of the graded limit.  The multiplicity of a particular $\lie g$-type is given by the number of certain lattice points in a convex polytope. Moreover, considering a particular face of that polytope encodes in fact the graded multiplicity.\\\\
A comparable study of HL-modules in other types is only partially explored. A first step was taken in \cite{CDM19} in type $D_n$ but it does not capture all the prime objects in the category $\cal F_\kappa$. There are important differences from the $A_n$ case and some new ideas seem to be necessary.
\subsubsection{Further Remarks} As we said, there are other subcategories of representations of $\cal F_q$ which were shown by Hernandez-Leclerc to be  monoidal categorifications of (infinite rank) cluster algebras. However, it is far from clear what subset of $\cal P^+_q$ is an index set for the prime representations corresponding to the cluster variables. Hence little is known about the characters or the graded limits of these representations.\\\\
Another example of prime representations comes from the theory of snake modules studied in types $A_n$ and $B_n$ in \cite{MY12, BM14}. Again the problem of studying the graded limits of these modules is wide open.\\\\

\section{Demazure modules,   Projective modules and Global Weyl modules}
Our focus in this section will be on the study of graded representations of $\lie g[t]$. We begin by establishing the correct category $\cal G$ of representations of the current algebra and introduce the projective objects and the global Weyl modules.  We then relate the study of local Weyl modules in Section \ref{classgrad}  to the $\lie g$-stable Demazure modules introduced in Section \ref{stable}. Next we  discuss the characters of the local Weyl modules and relate them to Macdonald polynomials. Finally, we discuss BGG-type reciprocity results. We conclude the section with some comments on the more recent work of \cite{FFD}, \cite{FM17a}, \cite{FD} \cite{Ka18}.
\subsection{The category $\cal G$} The study of this category was initiated in \cite{CG11} and we recall several ideas from that paper.  Recall from Section \ref{gradedlt} that we have a $\mathbb{Z}_+$-grading on $\lie g[t]$ and its universal enveloping algebra.
Define $\cal G$ to be the category whose objects are $\mathbb Z$-graded representations $V=\oplus_{m\in\mathbb Z} V[m]$ of $\lie g[t]$ with $\dim V[m]<\infty$ for all $m\in\mathbb Z$.  
The morphisms in the category are grade preserving maps of $\lie g[t]$-modules. \\\\
 Define  the restricted dual of an object $V$ in $\mathcal{G}$ by 
$$V^*=\bigoplus_{m\in\mathbb{Z}} V^{*}[m],\ \ \ V^{*}[m]=V[-m]^*.$$ Clearly $V^*$ is again an object of $\cal G$.\\\\
For any object $V$ of $\cal G$, each graded subspace $V[m]$ is a  finite-dimensional $\lie g$-module and we define the graded $\lie g$-character of $V$ to be the element of $\mathbb Z[P][[q^{\pm 1}]]$:
\begin{equation*}
\begin{split}
\ch_{\gr}V & =\sum_{\lambda\in P^+}\sum_{m\in\mathbb Z}\dim\Hom_{\lie g}(V(\lambda), V[m])q^m\ch V(\lambda)\\
& =\sum_{\mu\in P}\sum_{m\in\mathbb Z}\dim V[m]_\mu q^m e_{\mu}=\sum_{\mu\in P}p_{\mu}(q)e_\mu,\ \ p_\mu(q)\in\mathbb Z_+[[q^{\pm 1}]].
\end{split}
\end{equation*}

It is clear that for all $r\in\mathbb Z$ we have $\ch_{\gr}(\tau_r V)=q^r\ch_{\gr} V$, where $\tau_r$ is as defined in Section \ref{levelzero}).

\medskip
Finally, note that $\cal G$ is an abelian category and is closed under taking restricted duals. If $V$ and $V'$ are objects of $\cal G$ then $V\otimes V'$ is again an object in $\cal G$ if $\dim V<\infty$. 

\subsubsection{Finite-dimensional objects of $\cal G$}\label{findim} It is straightforward that if $V$ is a simple object of $\cal G$, then $V$ is concentrated in a single grade. In particular $V$ must be a finite-dimensional irreducible $\lie g$-module. In other words
$V\cong\tau_m\ev_0 V(\lambda)$ for some $\lambda\in P^+$ where $\ev_0$ is the evaluation $\lie g[t]\to \lie g$, $x\otimes t^r\mapsto\delta_{0,r}x$. From now we set $$V(\lambda,m)=\tau_m\ev_0 V(\lambda).$$
Another example of finite-dimensional modules in $\cal G$ are the $\lie g$-stable Demazure modules $V_w(\lambda)$, $\lambda\in\widehat P^+$ (see Section \ref{stable}) and the local Weyl modules  studied in Section \ref{locweyldef}.  We give a direct definition of those objects as $\lie g[t]$-modules here for the reader's convenience, and we also drop the $\tilde{}$ for ease of notation.
\\\\
Given $\lambda\in P^+$ the local Weyl module $W_{\loc}(\lambda)$  is the $\lie g[t]$-module generated by an element $v_\lambda$ and relations:$$x_i^+v_\lambda=0,\ \ (h\otimes t^r)v_\lambda=\delta_{r,0}\lambda(h)v_\lambda,\ \ (x_i^-)^{\lambda(h_i)+1}v_\lambda=0,\ \ \ i\in I,\ \ h\in\lie h.$$ Setting $\gr v_\lambda=r$ we see that $W_{\loc}(\lambda)$ can be regarded  as an object of $\cal G$ and we denote this as $W_{\loc}(\lambda, r)$. Clearly $W_{\loc}(\lambda, r)=\tau_rW_{\loc}(\lambda, 0)$. It was proved in \cite{CP01} that the local Weyl modules are finite-dimensional with unique irreducible quotients.\\\\
Given $\mu\in P^+$ and $\ell\in\mathbb N$, let $D(\ell,\mu)$ be the quotient of $W_{\loc}(\mu)$ by the submodule generated by elements $$(x_\alpha^-\otimes t^{s_\alpha-1})^{m_\alpha+1}v_\mu,\ \text{ if $m_{\alpha}<d_{\alpha}\ell$},\ \ (x_\alpha^-\otimes t^{s_\alpha})v_\mu,\ \ \alpha\in R^+$$
where $s_{\alpha}$ and $m_{\alpha}$ are determined by
$$ \mu(h_\alpha)=(s_\alpha-1)d_\alpha\ell+m_\alpha,\ \ 0<m_\alpha\le d_\alpha\ell.$$
The following was proved in \cite{CV13}.\begin{prop}\label{g-stable}
Suppose that $\lambda\in\widehat P^+$ and $w\in\widehat W$ is such that $w\lambda(h_i)\le 0$ for all $i\in I$ and assume that $\lambda(c)=\ell$. The module $V_w(\lambda)$ is isomorphic to $\tau_rD(\ell,\mu)$ where $\mu\in P^+$ is given by $\mu(h_i)=-w_\circ w\lambda(h_i)$ and $r=w \lambda(d)$.\hfill\qedsymbol
\end{prop}
\begin{rem}
An analogous presentation of non-stable Demazure modules is given in \cite{KV21} and we discuss this in the next section. These modules however are not objects of $\cal G$. \hfill\qedsymbol\end{rem}
\subsubsection{Relation between local Weyl and Demazure modules} The following corollary of Proposition \ref{findim} is easily established.
\begin{cor}\label{adelocdem}
If $\lie g$ is simply-laced then $$D(1,\mu)\cong W_{\loc}(\mu),\ \ \mu\in P^+.$$ 
\end{cor}
\begin{pf}
It suffices to prove that the following relation holds in $W_{\loc}(\mu)$: $$(x_\alpha^-\otimes t^{\mu(h_\alpha)})v_\mu=0,\ \ \alpha\in R^+.$$ But this follows by using $$(x_\alpha^-\otimes t^{\mu(h_\alpha)})v_\mu=(x_\alpha^+\otimes t)^{\mu(h_\alpha)}(x_\alpha^-)^{\mu(h_\alpha)+1}v_\mu=0.$$
Here the first equality is established by a simple calculation and using the relations in $W_{\loc}(\mu)$. The second equality holds since $W_{\loc}(\mu)$ is finite-dimensional and $$x_\alpha^+v_\mu=0\implies (x_\alpha^-)^{\mu(h_\alpha)+1}v_\mu =0.$$
\end{pf}
In the non-simply laced case it is not true in general that $W_{\loc}(\mu)\cong D(1,\mu)$. However it was proved in \cite{Na11} that $W_{\loc}(\mu)$ admits a decreasing filtration where the successive quotients are isomorphic to $\tau_rD(1,\mu_r)$ for some $r\in\mathbb Z$ and $\mu_r\in P^+$. In fact one can make a more precise statement which can be found in Section \ref{demflags}.
 \subsubsection{Projective modules and Global Weyl modules} Given $(\lambda,r)\in P^+\times\mathbb Z$ set $$P(\lambda, r)=\bu(\lie g[t])\otimes_{\mathbf{U}(\lie g)}V(\lambda, r).$$ It is not hard to check that $P(\lambda,r)$ is an indecomposable projective object of $\cal G$ and that there exists a surjective map $P(\lambda,r)\to V(\lambda, r)\to 0$ of $\lie g[t]$-modules. Equivalently $P(\lambda, r)$ is the $\lie g[t]$-module generated by an element $v_\lambda$ of grade $r$ subject to the relations:$$x_i^+v_\lambda=0,\ \ hv_\lambda=\lambda(h)v_\lambda,\ \ (x_i^-)^{\lambda(h_i)+1}v_\lambda=0,\ \ i\in I, h\in\lie h.$$
 The global Weyl module $W(\lambda, r)$ is 
the maximal quotient of $P(\lambda, r)$ such that $\mathrm{wt} W(\lambda,r)\subseteq\lambda-Q^+$. Equivalently it is the quotient of $P(\lambda,r)$ obtained by imposing the additional relations $(x_i^+\otimes t^k)v_\lambda=0$ for all $i\in I$ and  $k\ge 0$. Clearly we have the following sequence of surjective maps $$P(\lambda,r)\twoheadrightarrow W(\lambda ,r)\twoheadrightarrow W_{\loc}(\lambda,r)\twoheadrightarrow V(\lambda,r).$$
\subsubsection{The algebra $\mathbb A_\lambda$ and the bimodule structure on $W(\lambda, r)$}\label{alambda} Let $\lie h[t]_+=\lie h\otimes t\mathbb C[t]$ and for $\lambda\in P^+$ and $v_\lambda\in W(\lambda,r)_\lambda$ non-zero of grade $r$ let
$$\mathbb I_\lambda=\{u\in\bu(\lie h[t]_+): uv_\lambda=0\},\ \ \mathbb A_\lambda=\bu(\lie h[t]_+)/\mathbb I_\lambda.$$ Clearly $\mathbb A_\lambda$ is commutative and graded. 
Moreover $W(\lambda,r)$ is a $(\lie g[t],\mathbb A_\lambda)$-bimodule where the right action of $\mathbb A_\lambda$ is given by:$$(gv_\lambda)a=gav_\lambda,\ \ g\in\bu(\lie g[t]),\ \ a\in\mathbb A_\lambda.$$ 
To see that the action is well-defined, one must prove that
$$(\lie n^+\otimes \mathbb{C}[t])(h\otimes f)v_{\lambda}=0,\ \ (h'-\lambda(h'))(h\otimes f)v_{\lambda}=0,\ \ (x_i^-)^{\lambda(h_i)+1}(h\otimes f)v_\lambda=0$$
for all $i\in I$, $h,h'\in\lie h$ and $f\in\mathbb{C}[t]$. However, all relations are immediate to check. 
It was proved in \cite{CP01} (for the loop algebra; the proof is essentially the same for the current algebra) that $\mathbb A_\lambda$ can be realized as a ring of invariants as follows. Consider the polynomial ring $\mathbb C[x_{i,r}:i\in I, 1\le r\le\lambda(h_i)]$. The direct product of symmetric groups $$\mathcal S_\lambda=S_{\lambda(h_1)}\times\cdots\times S_{\lambda(h_n)}$$ acts on this ring in an obvious way and we have $$\mathbb A_\lambda\cong \mathbb C[x_{i,r}:i\in I, 1\le r\le\lambda(h_i)]^{\cal S_\lambda}.$$ The grading on  $\mathbb A_\lambda$ is given by requiring the grade of  $x_{i,r}$ being $r$.
Let $\bold I_\lambda$ be the maximal graded ideal in $\mathbb A_\lambda$.
The local Weyl module can then be realized as follows:
$$W_{\loc}(\lambda,r)= W(\lambda,r)\otimes_{\mathbb A_\lambda} {\mathbb A_\lambda}/{\bold I_\lambda}.$$
A nontrivial  consequence of the dimension conjecture discussed in Section \ref{locweyldef} (see \cite{CP01}, \cite{CFK10} for more details) is the following result.
\begin{prop}
  The global Weyl module $W(\lambda)$  is a free $\mathbb A_\lambda$-module of rank equal to the dimension of $W_{\loc}(\lambda)$.\hfill\qedsymbol
\end{prop}
The algebra $\mathbb A_\lambda$ plays an important role in the rest of the section.
\subsection{The category $\cal O$ for $\lie g$}\label{section415} Before continuing our study of the category $\cal G$, we discuss briefly, the resemblance of the theory with that of  the well-known category $\cal O$ for semi-simple Lie algebras. \\\\
The objects of  $\cal O$  are finitely generated weight  modules (with finite-dimensional weight spaces)  for $\lie g$ which are locally nilpotent for the action of $\lie n^+$. The morphisms are just $\lie g$-module maps. Given $\lambda\in \lie h^*$ one can associate to it
a Verma module $M(\lambda)$ which is defined as $$M(\lambda)=\bu(\lie g)\otimes_{\bu(\lie b)}\mathbb Cv_\lambda, $$ where $\mathbb Cv_\lambda$ is the one-dimensional $\lie b$-module given by $hv_\lambda=\lambda(h)v_\lambda$ and $ \lie n^+v_\lambda=0$. It is not hard to prove that $M(\lambda)$ is infinite-dimensional and has a unique irreducible quotient denoted by $V(\lambda)$ and any irreducible object in $\cal O$ is isomorphic to some $V(\lambda)$.  Moreover $V(\lambda)$ is finite-dimensional if and only if $\lambda\in P^+$; in particular $M(\lambda)$ is reducible if $\lambda\in P^+$. \\\\
The modules $M(\lambda)$ have  finite length and the multiplicity of $V(\mu)$ in the Jordan--Hölder series of $M(\lambda)$ is denoted by $[M(\lambda) : V(\mu)]$. The study of these multiplicities has been of great interest and there is extensive literature on the subject. Perhaps the starting point for this study is the famous result of Bernstein--Gelfand--Gelfand (BGG) which we now recall.\\\\
The category $\cal O$ has enough projectives, which means that for $\lambda\in\lie h^*$ there exists an indecomposable  module $P(\lambda)$ which is projective in $\cal O$ and we have surjective maps $$P(\lambda)\twoheadrightarrow M(\lambda)\twoheadrightarrow V(\lambda).$$  The following theorem (known as BGG-reciprocity) was proved in \cite{BGG}.
\begin{thm}
Given $\lambda_0\in\lie h^*$ there exist $\lambda_1,\dots,\lambda_r\in\lie h^*$ such that the  module $P(\lambda_0)$ has a decreasing filtration $P_0=P(\lambda_0)\supseteq P_1\supseteq P_2\supseteq \cdots \supseteq P_r\supseteq P_{r+1}= \{0\},$ and $$ P_i/P_{i+1}\cong M(\lambda_i),\ \ 0\le i\le r.$$ 
Moreover if we let $[P(\lambda): M(\mu)]$ be the multiplicity of $M(\mu)$ in this filtration then we have $[P(\lambda): M(\mu)]=[M(\mu): V(\lambda)].$\hfill\qedsymbol
\end{thm}

\begin{rem}
Although the filtration is not unique in general, a  comparison of formal characters shows that the filtration length and the multiplicity $[P(\lambda), M(\mu)]$ (see \cite[Section 3.7]{Hu08}) is independent of the choice of the filtration.
\end{rem}

More generally  a module in $\cal O$ which admits a decreasing sequence of submodules where the successive quotients are Verma modules is said to admit a standard filtration.
In the rest of this section we shall discuss an analog of this result for current algebras. \\\\
We will also explore other ideas stemming from the formal similarity between $\cal O$ and $\cal G$. For instance it is known that $\dim\Hom_{\lie g}(M(\lambda), M(\mu))\le 1$ and that any non-zero map between Verma modules is injective and we shall discuss its analog for current algebras. We shall also discuss an analog of tilting modules; in the category $\cal O$ these are defined to be modules which admit a filtration where the successive quotients are  Verma modules and also a filtration where the successive quotients are the restricted duals of Verma modules. It is known that  for each $\lambda\in\lie h^*$ there exists a unique indecomposable tilting module which contains a copy of $M(\lambda)$.
\subsection{BGG reciprocity in $\cal G$}\label{calbgg} In the category $\cal G$ the role of the Verma module is played by the global Weyl module. However, in general the global Weyl module $W(\lambda,r)$, $\lambda\in P^+$ does not have a unique finite-dimensional quotient in $\cal G$; for instance the modules $W_{\loc}(\lambda,r)$ and $V(\lambda,r)$ are usually not isomorphic and we have $$W(\lambda,r)\twoheadrightarrow W_{\loc}(\lambda,r)\twoheadrightarrow V(\lambda,r).$$ However both quotients have a uniqueness property; $ W_{\loc}(\lambda,r)$ is unique in the sense that  any finite-dimensional quotient of $W(\lambda,r)$ is actually a quotient of $W_{\loc}(\lambda,r)$ and $V(\lambda,r)$ is the  unique irreducible quotient of $W(\lambda, r)$. The further difference from the category $\cal O$ situation is that the global Weyl module is not of finite length. In spite of these differences, one is still able to formulate the appropriate version of BGG-reciprocity. Such a formulation was  first conjectured  in \cite{Bennet-Chari-Manning} and proved there for $\lie{sl}_2[t]$. The result was proved in complete generality in \cite{CI15} for twisted and untwisted current algebras; as usual the case  of $A_{2n}^{(2)}$ is much more difficult and one has to work with the hyperspecial current algebra. A key ingredient in the proof is to relate the character of the local Weyl module to specializations of (non)symmetric Macdonald polynomials (see Section \ref{macdonald} for a brief review). \\\\
The following is the main result of \cite{CI15}. 
\begin{thm}\label{mainci15}
Let $(\lambda, r)\in P^+\times \mathbb Z_+$. The module $P(\lambda,r)$ admits a decreasing series of submodules: $P_0=P(\lambda,r)\supseteq P_1\supseteq P_2\supseteq\cdots $ such that $$P_i/P_{i+1}\cong W(\mu_i,s_i), \ \ {\rm{for\ some}}\ \  (\mu_i,s_i)\in P^+\times \mathbb Z_+,$$ and $$[P(\lambda, r): W(\mu_i,s_i)] = [W_{\loc}(\mu_i,s_i): V(\lambda,r)].$$ 
\hfill\qedsymbol
\end{thm}
\subsubsection{Tilting modules} \iffalse Another important and interesting family of indecomposable objects in the category $\mathcal{O}$ are tilting modules characterized in terms of standard filtrations. A tilting module $M$ admits a standard filtration and a natural question about them is to ask for the filtration multiplicities. It turns out that there is a sort of duality between tilting modules and the projective modules introduced in Section~\ref{section415}.
\begin{thm}
Let $\lambda$ be an antidominant weight and $u\leq w$ in the Bruhat ordering of $W$. There is a unique indecomposable tilting module $T(\lambda)$ in $\mathcal{O}$ with integral weights satisfying $\mathrm{dim} T(\lambda)_{\lambda}=1$ and all weights $\mu$ of $T(\lambda)$ satisfy $\mu\leq \lambda$ such that 
$$[T(w\cdot \lambda): M(u\cdot \lambda)]=[P(w_{\circ}w\cdot \lambda: M(w_{\circ}u\cdot \lambda)]$$
where $w\cdot \lambda=w(\lambda+\rho)-\rho$ and $\rho$ is half the sum of positive roots.
\hfill\qed
\end{thm}\fi 
We discuss the construction of tilting modules  and some of their properties. These ideas were developed  in \cite{BB14a,BC12a,BC15a} and one works  in a suitable subcategory of $\cal G$. Thus, 
 let $\mathcal{G}_{\mathrm{bdd}}$ be the full subcategory of objects $M$ of $\mathcal{G}$ such that $M[j]=0$ for all $j\gg0$ and 
  $$\mathrm{wt}(M)\subseteq \bigcup_{i=1}^s \mathrm{conv} \ W \mu_i,\ \ \mu_1,\dots,\mu_s\in P^+$$
  where $\mathrm{conv}\ W \mu$ denotes the convex hull of the Weyl group orbit $W \mu$.
   \iffalse For example $V(\lambda,r)$ is an object in $\mathcal{G}_{\mathrm{bdd}}$ and exhaust all irreducible objects in $\mathcal{G}_{}$ and $\mathcal{G}_{\mathrm{bdd}}$ respectively. For the purpose of studying tilting modules, another important family of indecomposable modules in $\mathcal{G}_{\mathrm{bdd}}$ are the local Weyl modules $W_{\mathrm{loc}}(\lambda,r)$ introduced above. They will play the role of standard modules. Although costandard modules are usually defined to be the duals of standard modules, in this setting they are given by the dual global Weyl modules $W(\lambda,r)^{*}$. These modules are infinite-dimensional and give the third important family of indecomposable objects in $\mathcal{G}_{\mathrm{bdd}}$. 
\fi
An object $M$ in the category $\mathcal{G}_{\mathrm{bdd}}$ is called \textit{tilting} if it admits two increasing filtrations: \begin{gather*}M_0\subseteq M_1\subseteq\cdots  \subseteq M_r\subseteq \cdots,\  \ \  \ \ 
M^0\subseteq M^1\subseteq\cdots \subseteq M^r\subseteq\cdots \\\\ M=\bigcup_{r\ge 0} M_r= \bigcup_{r\ge 0}M^r,   \end{gather*}  such that $M_{i+1}/M_{i}$ (resp. $M^{i+1}/M^{i})$ is isomorphic to a finite direct sum of modules of the form $W_{\loc}(\lambda,r)$ (resp. to a sum of  dual global Weyl modules $W(\lambda, r)^*$) where $(\lambda,r)\in P^+\times\mathbb Z$.
One can also work with a dual definition of tilting modules, where one requires that the module has decreasing filtrations and the successive quotients are isomorphic to the dual local Weyl modules and the global Weyl modules respectively. 
\\\\
The following was proved in \cite[Section 2]{BC12a}.
\begin{thm} For $(\lambda,r)\in P^+\times \mathbb{Z}$ there exists an indecomposable  tilting module $T(\lambda,r)$ in  $\mathcal{G}_{\mathrm{bdd}}$  which maps onto the local Weyl module $W_{\loc}(\lambda,r)$ and such that $$\tau_r T(\lambda,0)= T(\lambda,r),\ \ T(\lambda,r)\cong T(\mu,s) \Leftrightarrow (\lambda,r)=(\mu,s).$$
Any indecomposable tilting module in $\mathcal{G}_{\mathrm{bdd}}$ is isomorphic to $T(\mu,s)$ for some $(\mu,s)\in P^+\times \mathbb{Z}$ and any tilting module in $\mathcal{G}_{\mathrm{bdd}}$ is isomorphic to a direct sum of indecomposable tilting modules. 
\hfill\qedsymbol
\end{thm}
The proof of the theorem relies on the following  necessary and sufficient condition for an object of  $\mathcal{G}_{\mathrm{bdd}}$
 to admit a filtration by dual global Weyl modules. Namely:\\\\
 $M$ admits a filtration by costandard modules if and only if  $\mathrm{Ext}^1_{\mathcal{G}}(W_{\mathrm{loc}}(\mu,s), M)=0,\ \forall (\mu,s)\in P^+\times \mathbb{Z}.$ 
 \begin{rem}
 This equivalence was established in \cite{BC12a} in the case when $\lie g$ is of type $A$. This is because the proof depended on knowing Theorem \ref{mainci15} which at that time had only been proved when $\lie g$ is of type $A$. However the proof given there goes through verbatim for any $\lie g$.\\\\
 \end{rem}
 \subsubsection{Tilting modules for $\lie{sl}_2$ and in Serre subcategories} The existence of tilting modules is proved in a very abstract way.  In \cite{BC15a} an explicit realization of the dual  modules  was given  in the case of $\lie{sl}_2$. In this case we identify $P^+$ with $\mathbb Z_+$. Recall also the algebra $\mathbb A_\lambda$ defined in Section \ref{alambda}; in this special case it is just the ring of symmetric polynomials in $\lambda$-variables.
 \begin{thm}
 Suppose that $\lie g=\lie{sl}_2$ and let $(\lambda, r)\in \mathbb Z_+\times\mathbb Z_+.$
The module $T(\lambda,r)^*$  is a free right $\mathbb{A}_{\lambda}$-module and $T(\lambda,r)\cong \tau_r \ T(\lambda,0)$. Moreover,
$$T(1,0)^*\cong W(1,0),\ \ T(\lambda,0)^*\cong \tau_{-r_{\lambda}} \ {\bigwedge}^{\lambda} W(1,0),\ \ r_{\lambda}=\binom{\lambda}{2},\ \ \lambda\geq 2.$$
$$ \mathrm{ch}_{\mathrm{gr}}T(\lambda,0)^*=\sum_{s=0}^{\lfloor\frac{\lambda}{2}\rfloor}t^{s(s-\lambda)}(1:t)_{s} \ch_{\mathrm{gr}}W(\lambda-2s,0)=t^{r_\lambda}(1:t)_{\lambda} \ch_{\mathrm{gr}}W_{\mathrm{loc}}(\lambda,r_{\lambda})^*$$
where
$$(1:t)_n=\frac{1}{(1-t)(1-t^2)\cdots (1-t^n)}$$

 \hfill\qedsymbol
\end{thm}
Very little is known about the structure or the character of the tilting modules in general. 
 \\\\
 A theory of tilting modules was also developed  for Serre subcategories of $\mathcal{G}$ which are defined as follows. Given a subset $\Gamma\subseteq P^+\times \mathbb{Z}$, we define a full subcategory  $\mathcal{G}(\Gamma)$ whose objects $M$ satisfy additionally 
$$[M: V(\lambda,r)]\neq 0 \Rightarrow (\lambda,r)\in \Gamma.$$ The category $\mathcal{G}(\Gamma)_{\mathrm{bdd}}$ is now defined in an obvious way. If $\Gamma=P^+\times J$, where $J$ is an (possibly infinite) interval in $\mathbb{Z}$, then the existence of tilting modules holds with $\mathcal{G}_{\mathrm{bdd}}$ and $P^+\times \mathbb{Z}$ replaced by  $\mathcal{G}(\Gamma)_{\mathrm{bdd}}$ and $\Gamma$ respectively (see \cite[Proposition 4.2 and Theorem 4.3]{BB14a}). The local and dual global Weyl modules in this setting are obtained by applying a certain natural functor to the standard and costandard modules in $\mathcal{G}$.

\subsubsection{Socle and Radical Filtration for local Weyl modules}
The local Weyl module $W_{\mathrm{loc}}(\lambda)$ has a natural increasing grading filtration induced from its graded module structure. This filtration coincides with the radical filtration (see \cite[Proposition 3.5]{KN12a}) which is defined as follows. For a module $M$ of $\mathbf{U}(\lie g[t])$ the radical filtration is given by
$$\cdots \subseteq \mathrm{rad}^k(M)\subseteq \cdots \subseteq \mathrm{rad}^1(M)\subseteq \mathrm{rad}^0(M)=M$$
where $\mathrm{rad}(M)$ is the smallest submodule of $M$ such that the quotient $M/\mathrm{rad}(M)$ is semi-simple and $\mathrm{rad}^k(M)$ is defined inductively by
$$\mathrm{rad}^k(M)=\mathrm{rad}(\mathrm{rad}^{k-1}(M)).$$
In particular, 
$$\mathrm{rad}^1(W_{\mathrm{loc}}(\lambda))=\bigoplus_{s>0}\mathbf{U}(\lie g[t])[s]v_{\lambda}.$$ There is another natural filtration on a module $M$, called the  socle filtration. It is  given as follows
$$0=\mathrm{soc}^0(M)\subseteq\mathrm{soc}^1(M)\subseteq \cdots\subseteq \mathrm{soc}^k(M)\subseteq \cdots,$$
where $\mathrm{soc}(M)=\mathrm{soc}^1(M)$ is the largest semi-simple submodule of $M$ and $\mathrm{soc}^k(M)$ is defined inductively by 
$$\mathrm{soc}^k(M)/\mathrm{soc}^{k-1}(M)=\mathrm{soc}(M/\mathrm{soc}^{k-1}(M)).$$
A module $M$ of $\mathbf{U}(\lie g[t])$ is called \textit{rigid} if the socle filtration coincides with the radical filtration. This is in particular the case, if $M$ is a finite-dimensional graded module such that $M/\mathrm{rad}(M)$ and $\mathrm{soc}(M)$ are both simple.  We remind the reader that when $\lie g$ is of type $ADE$ the local Weyl module is isomorphic to a level one Demazure module and hence embeds in a highest weight module for the affine Lie algebra. Given $\lambda\in P^+$, let $w\in\widehat W$ be such that $\Lambda=w^{-1}(w_0\lambda+\Lambda_0)\in\widehat P^+$. Since $V(\Lambda)$ is an irreducible integrable module for $\widehat{\lie g}$, it follows that any $\widehat{\lie b}$-submodule of $M=V_w(\Lambda)$  must contain the highest weight vector $v_\Lambda$. Hence any $\lie g[t]$-submodule of $M$ contains the $\lie g[t]$-module $U(\lie g[t])v_\Lambda$. In other words ${\rm{soc}}(M)$ must be simple and we must have ${\rm{soc}} (M)=U(\lie g[t])v_\Lambda\cong V(\Lambda|_{\lie h}, 0)$. So, we get
\begin{lem} Let $\lie g$ be of type $ADE$. Then ${\rm{soc} }(W_{\mathrm{loc}}(\lambda))\cong V(\Lambda|_{\lie h}, 0)$.
\hfill\qed
\end{lem}
It was proved in \cite{KN12a} that when $\lie g$ is of type  $ADE$ the local Weyl module is rigid. However, in general the socle of the local Weyl module is not simple and we give the  counterexample given in  \cite[Example 3.12]{KN12a}. In type $C_2$, the socle of the local Weyl module $W_{\mathrm{loc}}(2\omega_1+\omega_2)$ (the short root is $\alpha_1$) is isomorphic to $V(0, 3)\oplus V(\omega_2, 2)$.
\subsubsection{Maps between local Weyl modules} We apply the discussion on the socle of $W_{\loc}(\lambda,r)$ to study morphisms between local Weyl modules when $\lie g$ is of type $ADE$. For the purposes of this section it will be convenient to think of $W_{\loc}(\lambda,r)$ as modules for the subalgebra  $\lie g[t]\oplus\mathbb Cc\oplus\mathbb Cd$ of $\widehat{\lie g}$ where we let $c$ act as $1$ and the action of $d$ is given by the grading.  Recall the Bruhat order on $\widehat W$ given by $u\leq w$ if some substring of some reduced word for $w$ is a reduced word for $u$. 
%Let $\ell: \widehat W\to \mathbb Z_+$ be the length function. 
\begin{prop} Assume that $\lie g$ is of type $ADE$ and let $(\lambda,r), (\mu,s)\in P^+\times\mathbb Z$. Then,
$$\dim \mathrm{Hom}_{\mathcal{G}}(W_{\mathrm{loc}}(\lambda,r),W_{\mathrm{loc}}(\mu,s))\leq 1 $$ with equality holding if and only if there exist $w_1, w_2\in\widehat W$ and $\Lambda\in\widehat P^+$ such that the following hold:   $$ w_2\le w_1,\ \  w_1(\lambda+ \Lambda_0+r\delta)=\Lambda=w_2(\mu+\Lambda_0+s\delta).$$

Moreover, any non-zero map between local Weyl modules is injective.
\begin{proof} 
Let $\varphi: W_{\mathrm{loc}}(\lambda,r)\rightarrow W_{\mathrm{loc}}(\mu,s)$ be a non-zero homomorphism.  We first prove that this implies that  there exist $w_1, w_2\in\widehat W$ and $\Lambda\in\widehat P^+$ with $w_1(\lambda+\Lambda_0+r\delta)=\Lambda=w_2(\mu+ \Lambda_0+s\delta)$. To see this assume that $w_1(\lambda+\Lambda_0+r\delta)=\Lambda$ and $w_2(\mu+ \Lambda_0+s\delta)=\Lambda'$  with $\Lambda,\Lambda'\in \widehat{P}^+$ and let $ W_{\loc}(\mu,s)\hookrightarrow V(\Lambda')$ be the inclusion which exists since $ W_{\loc}(\mu,s)$ is isomorphic to a stable  Demazure module in $V(\Lambda')$. Since the image of $\varphi$ is non-zero it must include the simple socle of $W_{\loc}(\mu,s)$. This in turn implies that $\Lambda'$ is a weight of $W_{\loc}(\lambda, r)\hookrightarrow V(\Lambda)$. It follows that $\Lambda-\Lambda'$ must be a sum of affine positive roots. On the other hand since $\varphi(w_\lambda)\ne 0$ it follows that $\lambda+ \Lambda_0+r\delta$ must be a weight of $V(\Lambda')$ and hence $\Lambda$ is also a weight of $V(\Lambda')$. This forces $\Lambda'-\Lambda$ to be a sum of positive affine roots and also shows that $\Lambda=\Lambda'$. To see that $\varphi$ is injective, we note that otherwise both the kernel  and the cokernel of $\varphi$ would have to contain $v_\Lambda$ which is absurd.

Finally to see that the dimension of the homomorphism space is at most  one, it suffices to note that $\dim V(\Lambda)_{w\Lambda}=1$ for all $w\in\widehat W$.

\iffalse If we identify $W_{\mathrm{loc}}(\lambda)\cong V_w(\Lambda)$ and $W_{\mathrm{loc}}(\lambda)\cong V_u(\Lambda')$ we get a non-trivial induced map 
$$V_w(\Lambda)\rightarrow V(\Lambda')$$
which we also denote by $\varphi$ for convenience. In particular $\varphi(v_{\omega\Lambda})\neq 0$. Suppose that $w=s_i\tau$ and $\ell(w)>\ell(\tau)$ implying that $k:=\tau(\Lambda)(h_i)\geq 0$. Since $V_w(\Lambda)$ is $\lie g[t]$ stable and $\varphi$ commutes with the $\lie g[t]$ action we get
$$\varphi(v_{w\Lambda})=\varphi(x^{k}_{-\alpha_i}v_{\tau\Lambda})=x^{k}_{-\alpha_i}\varphi(v_{\tau\Lambda}).$$
Therefore, $\varphi(v_{\tau\Lambda})\neq 0$ and continuing in this way we get $\varphi(v_{\Lambda})\neq 0$. On the one hand this implies that $\varphi$ is injective (a non trivial kernel would contain the simple socle generated by $v_{\Lambda}$) and on the other hand we have $\Lambda=\Lambda'$ since $\varphi(v_{\Lambda})$ is a highest weight vector in $V(\Lambda')$ which is unique up to scalar. Thus $\varphi$ is also unique up to a scalar.\fi
\end{proof}
\end{prop}

\subsubsection{Morphisms between global Weyl modules} The study of the homomorphism space between global Weyl modules has also been studied in \cite{BCGM11} and confirms further the phenomenon that the global Weyl module plays a role similar to that of the Verma modules in category $\mathcal{O}$. The following result can be found in \cite[Theorem 3]{BCGM11}.
\begin{thm}
Let $\lambda, \mu\in P^+$ and  assume that   $\mu(h_i)=0$ for all $i\in I$ with $\omega_i(h_{\theta})\neq 1$. Then 
$$\mathrm{Hom}_{\mathcal{G}}(W(\lambda),W(\mu))=0,\ \ \text{if } \lambda\neq \mu,\ \ \mathrm{Hom}_{\mathcal{G}}(W(\mu),W(\mu))\cong \mathbb A_\mu.$$ Moreover any non-zero map $\varphi: W(\mu)\to W(\mu)$ is injective.\hfill\qed
\end{thm}
The restriction on $\mu$ is necessary (see \cite[Remark 6.1]{BCGM11}). For instance, in types $B_n$ and $D_n$ ($n\geq 6$)  we have  $\mathrm{Hom}_{\mathcal{G}}(W(\omega_2),W(\omega_4))\ne 0$.
 However the second statement, namely the injectivity of any non-zero map, is still expected to hold in general. 

\begin{rem} This theorem is quite unlike the analogous theorem  for local Weyl modules which was discussed in the preceding section.\end{rem}

  \subsection{Generalized Weyl modules,  Global Demazure modules and other directions.}

We now  discuss generalizations of some of the ideas presented earlier in this section. This is a brief and far from complete discussion of the papers of \cite{FM17a}, \cite{FD}, \cite{FFD} and we refer the interested readers to those papers for greater detail. We begin by elucidating the connection between local Weyl modules and specializations of Macdonald polynomials which was briefly mentioned in Section \ref{calbgg}. These polynomials are those associated with (anti) dominant weights.  We then discuss the work of \cite{FM17a} who introduced the notion of generalized Weyl modules for $\widehat{\lie n}^+$ and showed that their characters are again related to specializations of Macdonald polynomials associated with any integral weight.

We then move on to discuss the notion of global Demazure modules introduced by Dumanski and Feigin and state a few open problems regarding the homomorphism spaces between these objects. The aim is to generalize the global-local picture of Weyl modules for wider families of modules and  develop some modifications of results in this broader setting.
\subsubsection{Local Weyl modules, Generalized Weyl modules and Macdonald polynomials}\label{macdonald}   Let $$\mathcal{R}_{q,t}=\mathbb{Q}(q,t)[e_{\lambda}: \lambda\in P]\ \ \text{and }\  \mathcal{R}_q=\mathbb{Q}(q)[e_{\lambda}:\lambda\in P]$$ respectively be the group algebra of the weight lattice with coefficients in $\mathbb{Q}(q,t)$ and $\mathbb{Q}(q)$ respectively. Consider $R_{q,t}^{W}$ the subring of $W$-invariants where the action is induced from the action of $W$ on $P$ and define $R_{q}^{W}$ similarly. For $f\in R_{q,t}^{W}$, we denote by $[f]$ its constant term (i.e. the coefficient in front of $e_{0}$) and set
$$\displaystyle \nabla(q,t)=\prod_{\alpha\in (R \ + \ \mathbb{Z}_+\delta)} \frac{1-e_{\alpha}}{1-t^{-1}e_{\alpha}},\ \ e_{\delta}=q^{-1},\ \ \Delta(q,t)=\frac{\nabla(q,t)}{[\nabla(q,t)]}.$$ 
This ring $R_{q,t}^{W}$ and $R_q^{W}$ both admit a scalar product
$$\langle f,g\rangle_{q,t}=[f\overline{g}\Delta(q,t)]\ \ f,g\in R_{q,t}^{W},\ \ \ \langle f,g\rangle_{q}=[f\iota(g)\Delta(q,\infty)]\ \ f,g\in R_{q}^{W}$$
where $\overline{\cdot}$ is the involution on $R_{q,t}^{W}$ given by $t\mapsto t^{-1},\ q\mapsto q^{-1},\ e_{\lambda}\mapsto e_{-\lambda}$ and $\iota$ is the involution of $\mathcal{R}_q$ fixing $q$ and mapping $e_{\lambda}$ to $e_{-w_{\circ}\lambda}$.
Moreover, we have a natural basis $\{m_{\lambda}\}_{\lambda\in P^+}$ given by $$m_{\lambda}(q,t)=\sum_{\mu\in W\cdot \lambda} e_{\mu}$$
The \textit{symmetric Macdonald polynomials} $\{P_{\lambda}(q,t)\}_{\lambda\in P^+}$  are uniquely defined by the following two properties
\begin{enumerate}
    \item $\displaystyle P_{\lambda}(q,t)= m_{ \lambda}(q,t)+\sum_{\substack{\mu< \lambda\\ \mu\in P^+}} c_{\lambda,\mu} m_{\mu}(q,t),\ \ c_{\lambda,\mu}\in\mathbb{Q}(q,t),$ \vspace{0cm}
\item $\langle P_{\lambda}(q,t),P_{\mu}(q,t)\rangle_{q,t}=0,\ \ \lambda\neq \mu$
\end{enumerate}
These polynomials have the property that the limit $t\rightarrow \infty$ exists which we denote by $P_{\lambda}(q,\infty)=\lim_{t\rightarrow \infty} P_{\lambda}(q,t)$. The following result can be found in \cite[Theorem 4.2]{CI15}.  
\begin{thm} The family $\{P_{\lambda}(q,\infty)\}_{\lambda\in P^+}$ forms an orthogonal basis of $R^{W}_q$ with respect to the form $\langle \cdot,\cdot \rangle_q$.
Moreover
$$P_{\lambda}(q,\infty)=\ch_{\mathrm{gr}} W_{\loc}(\lambda),\ \ \lambda\in P^+$$\hfill\qed
\end{thm}
In the case when $\lie g$ is simply-laced it was already proved in \cite{I03} that the graded character of the stable Demazure module was given by the specialization of the Macdonald polynomial as in the above theorem; recall the connection between local Weyl modules and stable Demazure modules first made in \cite{CL06} in the case of  $\lie{sl}_{n+1}$ and then in \cite{FoL07} for $\lie g$ simply-laced.

At that time it was also known that this formula could not hold when $\lie g$ was not simply-laced. In the non-simply laced case the result of \cite{Na11} showed that the local Weyl module had a flag where the successive quotients were stable Demazure modules. The corresponding results for the twisted current algebras were studied in \cite{FK11}. However in the case of the twisted $A_{2n}^{(2)}$ one has to work with  a different current algebra \cite{CIK14}, called the hyperspecial current algebra.
\subsubsection{Nonsymmetric Macdonald polynomials} There is another family of polynomials $\{E_{\lambda}(q,t)\}_{\lambda\in P}$ indexed by the weight lattice called the \textit{nonsymmetric Macdonald polynomials}. They were introduced by Opdam \cite{Opdam95} and Cherednik \cite{C95}. First we define a new order on the set of weights $P.$ Consider the level one action of $\widehat{W}$ on $\lie h^*$ defined as follows (the action differs only for $s_0$)
 $$s_0\circ \mu:=s_{\theta}(\mu)+\theta,\ \mu\in \lie h^*.$$
Given $\lambda\in P$, we denote by $w_\lambda$ the unique minimal length element of $\widehat{W}$ such that $w_\lambda\circ \lambda$ is either miniscule or zero.
For $\lambda, \mu\in P$, we say $\mu <_b \lambda$ if and only if $w_\mu<w_\lambda$ with respect to the Bruhat order. 

Again we can define a scalar product $(\cdot,\cdot)_{q,t}$ on $\mathcal{R}_{q,t}$ and the family $\{E_{\lambda}(q,t)\}_{\lambda\in P}$ is  uniquely determined by the following two properties
\begin{enumerate}
\item 
$E_\lambda(q,t)=e_\lambda + \sum\limits_{\mu<_b\lambda}c_{\lambda, \mu} e_\mu, \, \, \, c_{\lambda,\mu}\in\mathbb{Q}(q,t)$, 
\item $(E_\lambda(q,t), e_\mu)_{q,t}=0 \ \ \text{if}\ \mu <_b \lambda.$
 \end{enumerate}
For a dominant weight $\lambda$ we have $P_{\lambda}(q,\infty)=\lim_{t\rightarrow \infty} E_{w_{\circ}\lambda}(q,t)= E_{w_{\circ}\lambda}(q,\infty)$ and hence by the above theorem the characters of local Weyl modules appear also as specializations of nonsymmetric Macdonald polynomials for anti-dominant weights
$$E_{w_{\circ}\lambda}(q,\infty)=\ch_{\mathrm{gr}} W_{\loc}(\lambda),\ \ \lambda\in P^+$$
 The natural question is whether other specializations are also meaningful in the sense that they have a representation theoretic interpretation. This leads to the definition of generalized local Weyl modules which can be found in \cite{FM17a}.\\

The reader should be warned that there are several notions of generalized Weyl modules, e.g. in \cite{CFK10,FKKS11,KL14} when the polynomial algebra is replaced by an arbitrary commutative algebra. But these are not the modules under consideration in this discussion.
\begin{defn}
Given $\mu\in P$ let $W_\mu$ be the $\widehat{\lie n}^+$-module generated by $v_\mu$ with relations, $$(h\otimes t^{r+1})v_\mu=0,\ r\ge0,\ \  \ (x^+_\alpha\otimes 1)^{\max\{-\mu(h_\alpha),0\}+1}v_\mu =0,\ \ (x^-_\alpha\otimes t)^{\max\{\mu(h_\alpha),0\}+1}v_\mu =0,\ \ \alpha\in R^+.$$
These are called the \textit{generalized local Weyl modules}.
\end{defn}
Note that for anti-dominant weights we obviously have $W_{\mu}\cong W_{\loc}(w_{\circ} \mu)$ as $\widehat{\lie n}^+$-modules and hence the character is again a specialized nonsymmetric Macdonald polynomial. The characters of $W_{\lambda}$ for $\lambda\in P^+$ are related to the Orr--Shimozono specialization of $E_{w_0\lambda}(q,t)$. The first part of the next proposition is proved in \cite{OS18} and the second part in \cite{FM17a}.
\begin{prop} Let $\lambda\in P^+$.
\begin{enumerate}
    \item The limit $E_{w_{\circ}\lambda}(q^{-1},0):=\lim_{t\to 0} E_{w_{\circ}\lambda}(q^{-1},t)$ exists and admits an explicit combinatorial formula in terms of quantum alcove paths.
    \item The character of $W_{\lambda}$ is given by $w_{\circ}E_{w_{\circ}\lambda}(q^{-1},0)$.
\end{enumerate}
\hfill\qed
\end{prop}

%\subsubsection{A construction of global Weyl modules and global Demazure modules}\label{construction}

\subsubsection{Recovering the global Weyl module from the local Weyl module} We recall a general construction which was first introduced in \cite{BCGM11}. Namely, if $V$ is any $\lie g[t]$-module, one can define an action of $\lie g[t]$ on $V[t]:=V\otimes\bc[t]$, by
$$(x\otimes t^r)(v\otimes t^s)=\sum_{j=0}^r \binom{r}{j}((x\otimes t^{r-j}) v)\otimes t^{s+j}.$$ In fact it was introduced in a more general context in Section 4 of \cite{BCGM11} by replacing $\mathbb C[t]$ by any commutative associative Hopf algebra $A$. Notice that if $V$ is generated by an element $v$ then $V[t]$ is generated by $v\otimes 1$. It was shown also that the fundamental global Weyl modules could be realized in this way by taking $V$ to be $W_{\loc}(\omega_i)$ for some $i\in I$. Moreover, it was shown in Proposition 6.2 of that paper that  if $\mu=\sum_{i=1}^ns_i\omega_i$ and $\mu(h_i)=0 $ if $\omega_i(h_\theta)\neq 1$ then
there exists an injective map $$W(\mu)\to W(\omega_1)^{\otimes s_1}\otimes\cdots\otimes W(\omega_n)^{\otimes s_n},$$ extending the assignment $w_\mu\to w_{\omega_1}^{s_1}\otimes\cdots\otimes  w_{\omega_n}^{s_n}$. This result was then established without the restriction on $\mu$ in \cite{Ka18} by different methods.

\subsubsection{A generalization} Suppose now that $W_1, \dots, W_r$ are graded $\lie g[t]$-modules with generators $w_1,\dots, w_r$. Then the associated global module is defined to be the submodule of $W_1[t]\otimes \cdots\otimes W_r[t]$ generated by the tensor product of $(w_1\otimes 1)\otimes\cdots \otimes (w_r\otimes 1)$. This notion was introduced in \cite[Section 1.3]{FM19g} and the resulting module is denoted as $R(W_1,\dots, W_r)$. In the case when the additional relations 
\begin{equation}\label{assa}(\lie h\otimes t\mathbb{C}[t])w_i = 0,\ \ hw_i=\lambda_i(h)w_i, \end{equation}
hold we can define (as in the case of global Weyl modules) a right action of $\mathbf{U}(\lie h[t])$ on $R(W_1,\dots,W_r)$. The algebra $\mathcal{A}(\lambda_1,\dots,\lambda_r)$ is as usual the quotient of $\mathbf{U}(\lie h[t]_+)$ by the annihilating ideal  of the cyclic vector $(w_1\otimes 1)\otimes\cdots \otimes (w_r\otimes 1)$. This algebra is harder to understand although one does have an  embedding
$$\mathcal{A}(\lambda_1,\dots,\lambda_r) \hookrightarrow \bigotimes_{i=1}^r \mathcal{A}(\lambda_i),\ \ \cal A(\lambda_i)\cong \mathbb{C}[z_i]$$
given by
$$(h\otimes t^m)\mapsto \sum_{i=1}^r (1\otimes \cdots\otimes ht^m\otimes \cdots \otimes 1)\mapsto \lambda_1(h)z_1^m+\cdots+ \lambda_r(h) z_r^m.$$  
Given $u\in\mathbb C$, set $$W_i(u)=W_i[t]\otimes_{\cal A(\lambda_i)}\mathbb C,$$ where we regard $\mathbb C$ as an $\cal A(\lambda_i)$-module by letting $z_i$ act as $u$.\\\\
The following was proved in \cite{FD}, \cite{FFD}.

\begin{thm} \label{thm2} Let $W_1,\dots ,W_r$ be as above. 
\begin{enumerit}
\item[(i)] The algebra $\mathcal{A}(\lambda_1,\dots,\lambda_r)$ is isomorphic to the subalgebra in $\mathbb{C}[z_1,\dots,z_r]$ generated by the polynomials $\lambda_1(h)z_1^m+\cdots+ \lambda_r(h) z_r^m,\ h\in\lie h,\ m\geq 0$.
\vspace{0,1cm}
\item[(ii)]  There exists a nonempty Zariski open subset $U$ of $\mathbb{C}^r$ (in particular $0\notin U$)  such that for all $\bold u=(u_1,\dots, u_r)\in U$
$$
R(W_1, \cdots , W_r) \otimes_{\mathcal{A}(\lambda_1, \dots , \lambda_r)} \mathbb{C}_{\mathbf{u}} \ \ \cong_{\lie g[t]} \ \ \bigotimes_{i = 1}^r W_i(u_i)
$$
where $\mathbb{C}_{\mathbf{u}}$ denotes the quotient of $\mathcal{A}(\lambda_1, \dots , \lambda_r)$ by the maximal ideal corresponding to $\mathbf{u}$.
\item[(iii)] The module $R(W_1, \dots , W_r)$ is a finitely-generated $\mathcal{A}(\lambda_1,\dots,\lambda_r)$-module.
\end{enumerit} \hfill\qedsymbol
\end{thm}
\iffalse 
\begin{proof} The first part is discussed before the theorem.  We consider the second part only for the$r=1$ case. In this situation we can choose $U=\mathbb{C}^{\times}$. So let $c\neq 0$ and consider the vector space isomorphism
$$W[t]\otimes_{\mathcal{A}(\lambda)} \otimes \mathbb{C}_c\simeq W(c),\ (v\otimes t^k)\otimes z\mapsto c^kz v$$
This is in fact an isomorphism of $\mathfrak{g}[t]$-modules which follows directly from \eqref{action1}:
\begin{align*}xt^m(v\otimes t^k\otimes 1)&=\sum_{j=0}^m \binom{m}{j} \big((xt^jv)\otimes t^{k+m-j}\big)\otimes 1&\\&=\sum_{j=0}^m  \binom{m}{j} \lambda(h)^{-1} \Big[(xt^jv)\otimes 1\Big]. (ht^{k+m-j})\otimes 1&\\&=\sum_{j=0}^m\binom{m}{j} c^{k+m-j}\big((xt^jv)\otimes 1\big) \otimes 1&\\& = c^k \Big[(x(t+c)^mv)\otimes 1\Big]\otimes 1\end{align*}

The algebras described in Theorem~\ref{thm2} (1) were studied in \cite{BCES16} and are referred to as \textit{algebras of deformed Newton sums}. One of their results is that $\bigotimes_{i=1}^r \mathcal{A}(\lambda_i)$ is finite over $\mathcal{A}(\lambda_1,\dots,\lambda_r)$ and that $\mathcal{A}(\lambda_1,\dots,\lambda_r)$ is noetherian. Since $\bigotimes_{i=1}^r W_i[t]$ is a finitely generated $\bigotimes_{i=1}^r \mathcal{A}(\lambda_i)$ module, it is also a finitely generated $\mathcal{A}(\lambda_1,\dots,\lambda_r)$ module. Hence also the submodule $R(W_1,\dots,W_r)$ which shows part (3).\end{proof}\fi
\subsubsection{}  The following conjecture of \cite{FD} generalizes the known results for local Weyl modules.
\begin{conj}\label{mainconj} The $\lie g[t]$-module
$$
R(W_1, \hdots , W_r) \otimes_{\mathcal{A}(\lambda_1, \dots , \lambda_r)} \mathbb{C}_0$$  is isomorphic to the fusion product of  the modules $W_i(c_i)$ 
for $(c_1,\dots,c_r)$ in some Zariski open subset of $\mathbb{C}^r$.
\end{conj}

In \cite{FD} the authors prove this conjecture for a certain families  of  Demazure modules when $\lie g$ is of type $ADE$, and in \cite{FFD} they drop the assumption on the type of $\lie g$. Given a collection of dominant integral weights $\underline{\lambda}=(\lambda_1,\dots,\lambda_r)$ we set $$\mathbb{D}_{\ell,\ell\underline{\lambda}}=R(D(\ell,\ell \lambda_1),\dots,D(\ell,\ell \lambda_r))$$
and let $v$ be the generating vector of $\mathbb D_{\ell,\ell\underline{\lambda}}$.  The following theorem has been proved for the tuple $\underline{\lambda}=(\omega_1,\dots,\omega_1,\dots,\omega_n,\dots,\omega_n)$ by Dumanski--Feigin and extended later by Dumanski--Feigin--Finkelberg to arbitrary tuples. \begin{thm} Let $\lambda\in P^+$ and $\underline{\lambda}=(\lambda_1,\cdots, \lambda_r)$ be such that $\lambda=\sum_{i=1}^r\lambda_i$. Then, we have an isomorphism 
$$D(\ell,\ell\lambda)\rightarrow \mathbb{D}_{\ell,\ell\underline{\lambda}}\otimes_{\mathcal{A}(\ell\lambda_1,\dots,\ell\lambda_r)}\mathbb{C}_0.$$
\hfill\qed
\end{thm}
An interesting question would be to determine the generators and relations for global Demazure modules. 

\begin{rem} Dumanski--Feigin--Finkelberg also prove that  $\mathbb{D}_{\ell,\ell\underline{\lambda}}$ is free over $\mathcal{A}(\ell\lambda_1,\dots,\ell\lambda_r)$ and that there exists a  tensor product decomposition 
$$\mathbb{D}_{\ell,\ell(\underline{\lambda}\cup \underline{\mu})}\otimes \mathbb{C}_{(c,d)}\cong (\mathbb{D}_{\ell,\ell \underline{\lambda}}\otimes \mathbb{C}_{c})\otimes (\mathbb{D}_{\ell,\ell \underline{\mu}}\otimes \mathbb{C}_{d})$$
provided that $c$ and $d$ have no common entries.  This is analogous to the well-known factorization of local Weyl modules which was proved in \cite{CP2001}.
\end{rem}
An interesting direction of research would be to study the homomorphisms between global Demazure modules and observe the analogues to homomorphisms between global Weyl modules discussed earlier in this section.

\subsection{} As we mentioned earlier the  current algebra is the derived algebra of the standard maximal parabolic subalgebra of the affine Lie algebra.  A natural problem is to  develop an analogous  theory for other parabolic subalgebras. This has been attacked for the first time in \cite{CKO18} and the two important families of local and global Weyl modules have been intensively studied, but many problems are still open. The global Weyl modules cotinue to be  parametrized by dominant integral weights of a semi-simple subalgebra of $\lie g$ depending on the choice of the maximal parabolic algebra.  However, the following interesting differences appear.
\begin{itemize}
    \item The algebra $\mathbb{A}_{\lambda}$ (modulo its Jacobson radical) is a Stanley--Reisner ring; in particular it has relations and is not a polynomial algebra (see \cite[Theorem 1]{CKO18}).
    \item The algebra $\mathbb{A}_{\lambda}$ and the global Weyl module can be finite-dimensional and this happens if and only if $\mathbb{A}_{\lambda}$ is a local ring.
    \item The global Weyl module is not a free $\mathbb{A}_{\lambda}$ module in general. However we expect the global Weyl module to be free over a suitable quotient algebra of $\mathbb{A}_{\lambda}$ corresponding to the coordinate ring of one of the irreducible subvarieties of $\mathbb{A}_{\lambda}$. 
\end{itemize}
The dimension of the local Weyl module depends on the choice of the maximal ideal of $\mathbb{A}_{\lambda}$. This was in the current algebra case one of the key observations together with the Quillen--Suslin theorem to obtain the freeness of global Weyl modules. It is still an open and interesting question to find the maximal ideals of $\mathbb{A}_{\lambda}$ producing the local Weyl modules of maximal dimension. An example has been discussed in \cite[Section 7.1]{CKO18}.

\section{Fusion product Decompositions, Demazure flags and  connections to combinatorics and Hypergeometric series}
In this section we collect together several results on Demazure modules which are of independent interest.
\subsection{Demazure modules revisited}
\subsubsection{A simplified presentation of Demazure modules.}
Recall that following \cite{Jos85, Ma88}, we gave in Section \ref{demazure}  a presentation of Demazure modules involving infinitely many relations. On the other hand we also discussed in Section \ref{adelocdem} that when $\lie g$ is simply-laced the local Weyl module $W_{\loc}(\mu,r)$ is isomorphic to a Demazure module occuring in a level one highest weight representation. The local Weyl module by definition has only finitely many relations. It turns out that this remains true for arbitrary Demazure modules. The following result was first proved in \cite{CV13} for $\lie g$-stable Demazure modules (see Proposition \ref{g-stable}) and was recently proved for arbitrary Demazure modules in \cite{KV21}.
\begin{thm}\label{mainthmgstable}
Suppose that $(\lambda, w)\in\widehat P^+\times \widehat W$ and assume that $\lambda(c)=\ell$, $\lambda(d)=r$  and $w\lambda|_{\mathfrak{h}} = \mu$. The module $V_w(\lambda)$ is isomorphic to
 a cyclic $\mathbf{U}(\widehat{\mathfrak{b}})$-module generated by a non-zero vector $v$ with the following relations: 
\begin{align*}
(h\otimes t^s)v=\delta_{s,0}\cdot \mu(h)v, \text{ for all } h\in \mathfrak{h},\, \,  dv = rv,\,\, cv = \ell v,
\end{align*}
and for $\alpha\in R^{\mp}(\mu)=\{\alpha\in R^+ : \mu(h_\alpha)\in \mp\mathbb{Z}_+\}$ we have
\begin{equation*}\label{demazurerelations1}
   \big(x_\alpha^{\pm}\otimes t^{s^{\pm}_\a-1}\big)^{m^{\pm}_\a+1}v=0, \ \text{if $m^{\pm}_{\alpha}<d_{\alpha}\ell$}; \ \ \big(x_\alpha^{\pm}\otimes t^{s^{\pm}_\a}\big)v=0,
  \end{equation*}
\begin{equation*}\label{demazurerelations2}
 \left(x_\alpha^+\otimes \mathbb{C}[t]\right)v=0,\ \left(x_\alpha^-\otimes t\right)^{\mathrm{max}\{0,\ \mu(h_{\alpha})-d_\a  \ell\}+1}v=0,\  \text{ if $\alpha\in R^+(\mu)$}
  \end{equation*}
\begin{equation*}\label{demazurerelations3}
\left(x_\alpha^-\otimes t\mathbb{C}[t]\right)v=0,\  \left(x_\alpha^+\otimes 1\right)^{-\mu(h_{\alpha})+1}v=0,\  \text{ if $\alpha\in R^-(\mu)$}
 \end{equation*}
where $s^{\pm}_\alpha,m^{\pm}_\alpha\in\bz_+$ are the unique integers such that  \begin{equation*}\label{ggg1}\mp\mu(h_{\alpha})=(s^{\pm}_{\alpha}-1)d_{\alpha}\ell+m^{\pm}_{\alpha},\ \ \ 0<m^{\pm}_{\alpha}\le d_{\alpha}\ell.\end{equation*} \qed
\end{thm}
\subsubsection{A tensor product theorem for $\lie g$-stable Demazure modules}
Recall that we discussed in Section \ref{locweyldef}
the realization of local Weyl modules as a fusion product of fundamental local Weyl modules. We also discussed in Section \ref{krtp}
the results of \cite{Na17} which gave the generators of the fusion products of Kirillov--Reshetikhin modules. These modules in the simply-laced case are known to be just $\lie g$-stable Demazure modules associated to weights of the form $\ell\lambda$ with $\ell\in\mathbb N$ and $\lambda\in P^+$. We remark here that in the simply-laced case, the results of \cite{Na17} are a vast generalization of the results of \cite{FoL07} where a presentation was given of  the fusion product of Demazure modules of a fixed level. This was achieved by showing that the fusion product was isomorphic to a Demazure module.  The following theorem 
which may be viewed as a Steinberg type decomposition theorem for $\lie g$-stable Demazure modules was proved in \cite{CSVW16} (see also \cite{VV21}) and completes the picture studied in \cite{FoL07}.
\begin{thm}\label{mainthmdecomposition}
  Let $\mathfrak{g}$ be a  finite-dimensional simple Lie algebra. Given $k\in \mathbb N$, let $\lambda\in P^+$, $\ell \in \mathbb{N}$ and suppose that $\lambda=\ell \ (\sum_{i=1}^k \lambda_i)+\lambda_{0}$ with $\lambda_0 \in P^+$ and $\lambda_i$ in the $\mathbb{Z}_+$-span of the $\omega_j^{\vee}$ for $1 \leq j \leq k$. Then there is an isomorphism of $\mathfrak{g}[t]$-modules
\begin{equation*}\label{eq:fact} D(\ell,\lambda)\cong  D(\ell,\lambda_0)^{z_0}*D(\ell,\ell\lambda_1)^{z_1}*\cdots*D(\ell,\ell\lambda_k)^{z_k}
\end{equation*}
where $z_0,\dots, z_k$ are distinct complex parameters. 
In particular the fusion product is independent of the choice of parameters. \hfill\qedsymbol
\end{thm}

The proof of the theorem relies on the simplified presentation of $D(\ell, \lambda)$ given in Theorem \ref{mainthmgstable} and a   character computation, using the Demazure character formula.  This allows one to  show that both modules in  the theorem  have the same $\lie h$-characters which is the crucial step to establish the theorem. The analogous theorem for  twisted  current algebras was proved in \cite{KV14}.

\subsubsection{} 
We discuss an interesting consequence of Theorem \ref{mainthmdecomposition}. Say that a  $\lie g[t]$-module is  {\em prime} if it is not isomorphic to a fusion product of non-trivial $\lie g[t]$-modules. The interested reader should compare this definition with that given in Section \ref{prime} where an analogous definition was made in the context of quantum affine algebras.
The following factorization result is a consequence of Theorem \ref{mainthmdecomposition}: 
\begin{cor} Given $\ell \geq 1$ and $\lambda \in P^+$ write  $$\lambda = \ell \, \left(\sum_{i=1}^n d_im_i \omega_i\right) + \lambda_0,\ \  
\lambda_0 = \sum_{i=1}^n r_i \omega_i,  \ \ 0 \leq r_i <d_i\ell, \ \ m_i\in\mathbb Z_+.$$
Then $D(\ell,\lambda)$ has the following fusion product factorization:    \begin{equation}\label{eq:primefact} 
  D(\ell,\lambda)\cong_{\lie g[t]}  D(\ell,\ell d_1\omega_1)^{*m_1}*D(\ell,\ell d_2\omega_2)^{*m_2}*\cdots*D(\ell,\ell d_n\omega_n)^{*m_n}* D(\ell,\lambda_0).
\end{equation}
In addition, if we assume that $\lie g$ is simply-laced then \ref{eq:primefact} gives prime factorization of $D(\ell,\lambda)$   (i.e., each module on the ritht hand side is prime).

\end{cor}

\iffalse \subsubsection{ Generalzied $Q$-systems} 
 The notion of  $Q$-systems introduced in \cite{HKOTY99}. A $Q$-system is a short exact sequence of $\lie g$-modules:
\[ 0 \rightarrow \bigotimes_{j\sim i} V(\ell \omega_j) \rightarrow V(\ell \omega_i) \otimes V(\ell \omega_i) \rightarrow  V((\ell+1) \omega_i) \otimes V((\ell-1) \omega_i) \rightarrow 0\]
where $\omega_i$ is a miniscule weight of $\lie g$ and $j\sim i$ means that nodes $i$ and $j$ are connected by an edge in the Dynkin diagram. Generalizations of $Q$-systems considered 
in \cite{CSVW16, CV13, FHD14, KV14} involve replacing the tensor products above by fusion products of certain $\mathring{\lie g}[t]$-modules.
The following generalization of $Q$-system was established in \cite[Section 5]{CSVW16} (see also \cite{VV21}):
\begin{prop}\label{prop:csvw-qsys} 
  Let $\lie g$ be simply-laced. Let $\ell \ge 1, \lambda \in P^+$ with $\ell \ge \lambda(h_\theta)$ and suppose $\omega_i$ is a miniscule weight such that $\lambda(h_i) >0$. Set $\mu = \ell \omega_i + \lambda - \lambda(h_i) \alpha_i$. Then there exists a short exact sequence of $\lie g[t]$-modules:
\begin{gather*} 0\to\tau_{\lambda(h_i)}\left( D(\ell, \ell\mu_1)*D(\ell, \mu_0)\right)\to D(\ell, \ell\omega_i)*D(\ell,\lambda) \to D(\ell+1,(\ell+1)\omega_i)*D(\ell+1 ,\lambda-\omega_i)\to 0.\end{gather*}
\end{prop}
\fi
\subsection{Demazure flags}\label{demflags}
In this section we explain a connection  between modules admitting Demazure flags and combinatorics and  hypergeometric series.\\\\
Say that a finite-dimensional $\lie g[t]$-module $M$ has a level $m$-Demazure flag if it admits a decreasing family of submodules, $M=M_0\supseteq M_1\supseteq\cdots \supseteq M_r\supseteq 0,$ such that $$M_j/M_{j+1}\cong\tau_{r_j}D(m, \mu_j),\ \ r_j\in\mathbb Z,\ \ \mu_j\in P^+.$$ It is not hard to see by working with graded characters, that if $M=M_0'\supseteq M_1'\supseteq\cdots \supset M_s'\supseteq 0$ is another level $m$-Demazure flag then $r=s$ and the multiplicity of $\tau_rD(m,\mu)$ in both flags is the same. Hence we define $[M: \tau_r D(m,\mu)]$ to be the number of times  $\tau_r D(m,\mu)$ occurs in a level $m$-Demazure flag of $M$.\\\\
 The study of Demazure flags goes back to the work of  Naoi \cite{Na11} on local Weyl modules in the non-simply laced case. It was proved in that paper that these modules admit a level one Demazure flag. This was done by first showing that in the simply-laced case every $\lie g[t]$-stable Demazure module of level $\ell$ admits  a Demazure flag of level $m$
   if $m \geq \ell \geq 1$. In the case when $\ell=1$ and $m=2$ the multiplicities occurring in this flag can be explicitly related to the multiplicity of the level one flag of the local Weyl module for non-simply laced Lie algebras.
   However, the methods do not lead to precise formulae for the multiplicity. In this section we discuss how one might approach this problem using different kinds of  generating series.

\subsubsection{The case of $\lie{sl}_2$ and level two flags} A first step to calculate the multiplicity was taken in \cite{CSSW14} for the Lie algebra $\mathfrak{sl}_2$ when $m=2$ and $\ell=1$. Then the graded multiplicities can be expressed by $q$-binomial coefficients \cite[Theorem 3.3]{CSSW14}:
\begin{equation*}
\begin{split}
[D(1,\mu\omega):D(2,\nu\omega)]_q & =\sum_{p\geq 0}[D(1,\mu\omega):\tau_p D(2,\nu\omega)]\ q^p\\ & =
\begin{cases} q^{(\mu-\nu)/2\lceil{\mu/2\rceil}}\qbinom{\lfloor\mu/2\rfloor}{(\mu-\nu)/2}_q,\ \ \mu-\nu\in 2\mathbb{Z}_+,\\ 0\ \ {\rm{otherwise}}.\end{cases}
\end{split}
\end{equation*}

\subsubsection{The case of $\lie{sl}_2$ and arbitrary level} A more general approach in the $\mathfrak{sl}_2$ case was taken in the articles \cite{BCSV15,BKu17} and in the $A_2^{(2)}$ case in \cite{BCK18a}.  In those papers, the  authors found a connection to algebraic combinatorics and number theory. We first need to introduce more notation. Define a family of  generating series  by
$$A^{\ell\rightarrow m}_{n}(x,q)=\sum_{k\geq 0}\ [D(\ell,(n+2k)\omega):D(m,n\omega)]_q \cdot x^k,\, \, n\ge 0.$$
Introduce the partial theta function $\theta(q,z)=\sum_{k=0}^{\infty} q^{k^2}z^k$ and let $$P_n(x) = \sum_{k =0}^{\floor{\frac{n}{2}}} (-1)^k \,{n-k \choose k} \, x^k$$ and note that the polynomials $P_n(x)$ are related to the Chebyshev polynomials $U_n(x)$ of the second kind as follows $P_n(x^2) = x^n \, U_n(\, (2x)^{-1}).$ The following theorem can be found in \cite[Theorem 1.6]{BCSV15} and \cite[Corollary 1.3]{BCSV15} respectively. 
\begin{thm}
\begin{enumerit}
\item[(i)]  Let   $n,m \in\bz_+$ and write  $n = ms + r$ with $s\in\mathbb Z_+$ and $0 \leq r < m$. Then 
$$A^{1\rightarrow m}_{n}(x,1)=\frac{P_{m-r-1}(x)}{P_m(x)^{s+1}}.$$ 

\item[(ii)] The specializations
$$A^{1\rightarrow 3}_{1}(1,q),\ \ \ q\cdot A^{1\rightarrow 3}_{1}(q,q),\ \ \ A^{1\rightarrow 3}_{0}(1,q^2)+q A^{1\rightarrow 3}_{2}(1,q^2),\ \ \ q^4A^{1\rightarrow 3}_{2}(q^2,q^2)+q A^{1\rightarrow 3}_{0}(q^2,q^2)$$
coincide with fifth order mock theta functions of Ramanujan.
\item[(iii)] The series $A^{1\rightarrow 2}_{n}(x,q)$ and $A^{2\rightarrow 3}_{n}(x,q)$ can be expressed as a linear combination of specializations of the partial theta function $\theta$ whose coefficients are given by products of $q$-binomial coefficients. 
\end{enumerit}
\hfill\qed
\end{thm}
\subsubsection{Combinatorics of Dyck paths and the functions $A_n^{1\to m}$} We further discuss the $\mathfrak{sl}_2$ case and its connection to the combinatorics of Dyck paths. In \cite{BKu17} a combinatorial formula has been obtained whose ingredients we will now explain. A \textit{Dyck path} is a diagonal lattice path from the origin $(0,0)$ to $(s,n)$ for some non-negative integrs $s,n\in \bz_+$,
such that the path never goes below the x-axis. We encode such a path by a 01-word, where $1$ encodes the up-steps and $0$ the down-steps.  We denote by $\mathcal{D}^N_{n}$ the set of Dyck paths that end at height $n$ and which do not cross the line $y=N$. The following picture is an example of an element in $\mathcal{D}^4_{1}$. 
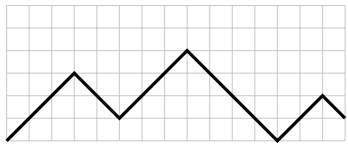
$$
\begin{tikzpicture}[scale=0.3]
\draw [step=1,thin,gray!40] (0,0) grid (15,6);
\draw [very thick] (0,0) -- (3,3) -- (5,1) -- (8,4)--(12,0)--(14,2)--(15,1);
\end{tikzpicture}$$
For $n,m\in \bz_+$ let $n_0,n_1\in\bz_+$ be such that $n_0<m$ and $n=m n_1+n_0$. If $n<m$ we set $A(m,n)=\emptyset$ and otherwise define
$$A(m,n):=\{(i_1,m),(i_2,m+1),\dots,(i_{n-m+1},n)\}\subseteq \bz_+^2$$
where $i_1<\cdots<i_{n-m+1}$ is the natural ordering of the set 
$$\{0,\dots,n\}\backslash \{pn_1+n_0+\min\{0,(p-1)-n_0\},\ 1\leq p\leq m\}.$$
Given a pair of non-negative integers $(a,b)\in \bz^2_+$, we say that $P\in \mathcal{D}_n^{\max\{m-1,n\}}$ is $(a,b)$-admissible if and only if $P$ satisfies the following property. If $P$ has a peak at height $b$, the subsequent path is strictly above the line $y=a$. For example, the path above is not $(0,3)$-admissible. \\\\
 Let $\mathcal{D}_{m,n}$ be the set  Dyck paths in $\mathcal{D}_n^{\max\{m-1,n\}}$, which are  $(a,b)$-admissible for all $(a,b)\in A(m,n)$. \\\\ 
The major statistics of a Dyck path was studied first by MacMahon \cite{MacMahon} in his interpretation of the $q$-Catalan numbers. Let $P=a_1\cdots a_{s}$, $a_i\in\{0,1\}$ be  a Dyck path of length $s$. The major and comajor index are defined by 
$$\text{maj}(P)=\sum_{\substack{1\leq i < s,\\ a_i>a_{i+1}}} i,\ \ \ \text{comaj}(P)=\sum_{\substack{1\leq i < s,\\ a_i>a_{i+1}}} (s-i).$$
The following was proved  in \cite[Theorem 4]{BKu17}.
\begin{thm}\label{dyck1} Let $m\in\bn$, $n\in\bz_+$. We have, 
$$A^{1\rightarrow m}_n(x,q)=\sum_{P\in \mathcal{D}_{m,n}}q^{\text{comaj}(P)}\ x^{d(P)}$$
where $d(P)$ denotes the number of down-steps of $P$.\hfill\qedsymbol
\end{thm}
In the twisted case graded and weighted generating functions encode again the multiplicity of a given Demazure module. For small ranks these generating functions are completely determined in \cite{BCK18a} and they define hypergeometric series and are related to the $q$–Fibonacci polynomials defined by Carlitz. For more details we refer the reader to \cite[Section 2]{BCK18a}.
\subsubsection{The general case} It is still an open problem to come up with closed or even recursive formulas for the generating series for other finite-dimensional simply-laced Lie algebras; the multiplicities and generating functions are defined in the obvious way. However, some progress has been made in \cite{BCSW21} for the Lie algebra $\mathfrak{sl}_{n+1}$ and the connection to Macdonald polynomials was established. The following result can be derived from \cite{BCSW21}.
\begin{thm} Let $\lie g=\mathfrak{sl}_{n+1}$ and $\lambda,\mu \in P^+$ such that $\lambda-\mu=\sum_{i=1}^n k_i \alpha_i, \ k_i\in \mathbb{Z}_+$. Then,  
$$[D(1,\lambda): D(2,\mu)]_q= \prod_{i=1}^n \big[D(1,(\mu(h_i)+2k_i)\omega): D(2,\mu(h_i)\omega)\big]_q$$
where $\omega$ is the corresponding fundamental weight for $\mathfrak{sl}_2$.
\hfill\qed
\end{thm}
So combining the above theorem with the combinatorial formula in Theorem~\ref{dyck1} gives a combinatorial formula for graded multiplicities of level 2 Demazure modules in level one flags.

\bibliographystyle{plain}
\bibliography{bibfile}

\end{document}